\documentclass{amsart}

\usepackage{amsfonts,amssymb,amsmath,amscd,amsthm}
\usepackage[mathscr]{eucal}
\usepackage{mathrsfs}
\usepackage{mathbbol}

\usepackage{oldgerm,units,eepic}
\usepackage{wrapfig,epsfig}

\usepackage{ifthen}

\newtheorem{theorem}{Theorem}[section]
\newtheorem{proposition}[theorem]{Proposition}
\newtheorem{definition}[theorem]{Definition}

\newtheorem{lemma}[theorem]{Lemma}

\newtheorem{corollary}[theorem]{Corollary}

\newtheorem{example}[theorem]{Example}
\newtheorem{remark}[theorem]{Remark}

\newtheorem{observation}[theorem]{Observation}

\newtheorem{axiom}[theorem]{Axiom}

\newcommand{\rTheorem}[2]{\thmSpc \noindent \emph{\textbf{Theorem #1:} #2}\thmSpc}

\newcommand{\Ref}[1]{(\ref{#1})}

\newcommand{\Comp}{\mathbb C}
\newcommand{\Real}{\mathbb R}
\newcommand{\Rati}{\mathbb Q}
\newcommand{\Int}{\mathbb Z}
\newcommand{\Net}{\mathbb N}

\newcommand{\Fld}{\mathbb K}

\newcommand{\one}{\mathbb{1}}
\newcommand{\zero}{\mathbb{0}}

\newcommand{\Trop}{\mathbb T}

\newcommand{\tUnit}{\Real^\nu}
\newcommand{\etUnit}{\bar{\Real}^\nu}
\newcommand{\eReal}{\bar{\Real}}

\newcommand{\trop}[1]{\mathcal{#1}}

\newcommand{\tA}{\trop{A}}
\newcommand{\tB}{\trop{B}}

\newcommand{\tI}{\trop{I}}

\newcommand{\tS}{\trop{S}}

\newcommand{\tZ}{\trop{Z}}

\newcommand{\tpF}{f}

\newcommand{\tfF}{\tilde{f}}

\newcommand{\bbx}{\square}
\newcommand{\To}{\longrightarrow }

\newcommand{\tUniS}{-\infty}

\newcommand{\uuu}[1]{#1^\nu}

\newcommand{\maxPlusToU}{\theta}

\newcommand{\homToU}{\nu}

\newcommand{\epiToMaxPlus}{\pi}

\newcommand{\invA}[1]{#1^{\bigtriangledown}}

\newcommand{\sinvA}[1]{#1^{\diamond}}

\newcommand{\io}{\iota}

\newcommand{\al}{\alpha}

\newcommand{\lm}{\lambda}
\newcommand{\Lm}{\Lambda}
\newcommand{\sig}{\sigma}

\newcommand{\itA}{\invA{\tA}}
\newcommand{\itB}{\invA{\tB}}

\newcommand{\stA}{\sinvA{\tA}}

\newcommand{\ptA}{\tA}

\newcommand{\Genr}[1]{\widetilde{#1}}
\newcommand{\gtI}{\Genr{\tI}}

\newcommand{\dt}{{_\bullet}}

\newcommand{\Val}{Val}
\newcommand{\nVal}{Val}

\newcommand{\TrS}{\oplus}
\newcommand{\TrP}{\odot}
\newcommand{\iTrP}{\TrP^{\triangledown}}

\newcommand{\tSR}{(\Trop,\TrS,\TrP)}

\newcommand{\eMaxPlus}{(\eReal,\max,+)}

\newcommand{\OP}{\left(}
\newcommand{\CP}{\right)}

\def \aaa{a}
\def \aab{b}
\def \aac{c}

\newcommand{\inva}[1]{#1^{\triangledown}}
\newcommand{\pa}{a}
\newcommand{\ia}{\inva{\aaa}}

\newcommand{\sia}{\sinvA{a}}

    \ifx\proof\undefined
    \newenvironment{proof}{
    \smallskip
    \noindent\emph{Proof.}}{\hfill\(\Box\)
    \bigskip
    } \fi

\newcommand{\vMat}[4]{
\OP \begin{array}{cc}
  #1 & #2 \\
  #3 & #4
\end{array}\CP}

\newcommand{\vvMat}[9]{\OP \begin{array}{ccc}
  #1 & #2 & #3\\
  #4 & #5 & #6\\
  #7 & #8 & #9\\
\end{array}\CP}

\newcommand{\bfem}[1]{\textbf{\emph{#1}}}

\newcommand{\ifdef}[3]{\ifthenelse{\equal{#1}{true}}{#2}{#3}}

\newcommand {\parSpc} {\vskip 0.1cm}
\newcommand {\thmSpc} {\vskip 0.1cm}

\def\trn{{\operatorname{t}}}
\def\idem{{\operatorname{idm}}}

\def\maps{F}

\def\bfi{ \textbf{i}}

\newcommand{\Adj}[1]{Adj({#1})}

\def\tpF{f}

\makeatletter

\def\diagram{\m@th\leftwidth=\z@ \rightwidth=\z@ \topheight=\z@
\botheight=\z@ \setbox\@picbox\hbox\bgroup}

\def\enddiagram{\egroup\wd\@picbox\rightwidth\unitlength
\ht\@picbox\topheight\unitlength \dp\@picbox\botheight\unitlength
\hskip\leftwidth\unitlength\box\@picbox}

\def\bfig{\begin{diagram}}
\def\efig{\end{diagram}}
\newcount\wideness \newcount\leftwidth \newcount\rightwidth
\newcount\highness \newcount\topheight \newcount\botheight

\def\ratchet#1#2{\ifnum#1<#2 \global #1=#2 \fi}

\def\putbox(#1,#2)#3{%
\horsize{\wideness}{#3} \divide\wideness by 2
{\advance\wideness by #1 \ratchet{\rightwidth}{\wideness}}
{\advance\wideness by -#1 \ratchet{\leftwidth}{\wideness}}
\vertsize{\highness}{#3} \divide\highness by 2
{\advance\highness by #2 \ratchet{\topheight}{\highness}}
{\advance\highness by -#2 \ratchet{\botheight}{\highness}}
\put(#1,#2){\makebox(0,0){$#3$}}}

\def\putlbox(#1,#2)#3{%
\horsize{\wideness}{#3}
{\advance\wideness by #1 \ratchet{\rightwidth}{\wideness}}
{\ratchet{\leftwidth}{-#1}}
\vertsize{\highness}{#3} \divide\highness by 2
{\advance\highness by #2 \ratchet{\topheight}{\highness}}
{\advance\highness by -#2 \ratchet{\botheight}{\highness}}
\put(#1,#2){\makebox(0,0)[l]{$#3$}}}

\def\putrbox(#1,#2)#3{%
\horsize{\wideness}{#3}
{\ratchet{\rightwidth}{#1}}
{\advance\wideness by -#1 \ratchet{\leftwidth}{\wideness}}
\vertsize{\highness}{#3} \divide\highness by 2
{\advance\highness by #2 \ratchet{\topheight}{\highness}}
{\advance\highness by -#2 \ratchet{\botheight}{\highness}}
\put(#1,#2){\makebox(0,0)[r]{$#3$}}}

\def\adjust[#1]{} 

\newcount \coefa
\newcount \coefb
\newcount \coefc
\newcount\tempcounta
\newcount\tempcountb
\newcount\tempcountc
\newcount\tempcountd
\newcount\xext
\newcount\yext
\newcount\xoff
\newcount\yoff
\newcount\gap%
\newcount\arrowtypea
\newcount\arrowtypeb
\newcount\arrowtypec
\newcount\arrowtyped
\newcount\arrowtypee
\newcount\height
\newcount\width
\newcount\xpos
\newcount\ypos
\newcount\run
\newcount\rise
\newcount\arrowlength
\newcount\halflength
\newcount\arrowtype
\newdimen\tempdimen
\newdimen\xlen
\newdimen\ylen
\newsavebox{\tempboxa}%
\newsavebox{\tempboxb}%
\newsavebox{\tempboxc}%

\newdimen\w@dth

\def\setw@dth#1#2{\setbox\z@\hbox{\m@th$#1$}\w@dth=\wd\z@
\setbox\@ne\hbox{\m@th$#2$}\ifnum\w@dth<\wd\@ne \w@dth=\wd\@ne \fi
\advance\w@dth by 1.2em}

\def\t@^#1_#2{\allowbreak\def\n@one{#1}\def\n@two{#2}\mathrel
{\setw@dth{#1}{#2}
\mathop{\hbox to \w@dth{\rightarrowfill}}\limits
\ifx\n@one\empty\else ^{\box\z@}\fi
\ifx\n@two\empty\else _{\box\@ne}\fi}}
\def\t@@^#1{\@ifnextchar_{\t@^{#1}}{\t@^{#1}_{}}}
\def\to{\@ifnextchar^{\t@@}{\t@@^{}}}

\def\t@left^#1_#2{\def\n@one{#1}\def\n@two{#2}\mathrel{\setw@dth{#1}{#2}
\mathop{\hbox to \w@dth{\leftarrowfill}}\limits
\ifx\n@one\empty\else ^{\box\z@}\fi
\ifx\n@two\empty\else _{\box\@ne}\fi}}
\def\t@@left^#1{\@ifnextchar_{\t@left^{#1}}{\t@left^{#1}_{}}}
\def\toleft{\@ifnextchar^{\t@@left}{\t@@left^{}}}

\def\two@^#1_#2{\allowbreak
\def\n@one{#1}\def\n@two{#2}\mathrel{\setw@dth{#1}{#2}
\mathop{\vcenter{\lineskip\z@\baselineskip\z@
                 \hbox to \w@dth{\rightarrowfill}%
                 \hbox to \w@dth{\rightarrowfill}}%
       }\limits
\ifx\n@one\empty\else ^{\box\z@}\fi
\ifx\n@two\empty\else _{\box\@ne}\fi}}
\def\tw@@^#1{\@ifnextchar _{\two@^{#1}}{\two@^{#1}_{}}}
\def\two{\@ifnextchar ^{\tw@@}{\tw@@^{}}}

\def\tofr@^#1_#2{\def\n@one{#1}\def\n@two{#2}\mathrel{\setw@dth{#1}{#2}
\mathop{\vcenter{\hbox to \w@dth{\rightarrowfill}\kern-1.7ex
                 \hbox to \w@dth{\leftarrowfill}}%
       }\limits
\ifx\n@one\empty\else ^{\box\z@}\fi
\ifx\n@two\empty\else _{\box\@ne}\fi}}
\def\t@fr@^#1{\@ifnextchar_ {\tofr@^{#1}}{\tofr@^{#1}_{}}}
\def\tofro{\@ifnextchar^ {\t@fr@}{\t@fr@^{}}}

\def\mon{\mathop{\m@th\hbox to
      14.6\P@{\lasyb\char'51\hskip-2.1\P@$\arrext$\hss
$\mathord\rightarrow$}}\limits} 
\def\leftmono{\mathrel{\m@th\hbox to
14.6\P@{$\mathord\leftarrow$\hss$\arrext$\hskip-2.1\P@\lasyb\char'50%
}}\limits} 
\mathchardef\arrext="0200       

\def\settypes(#1,#2,#3){\arrowtypea#1 \arrowtypeb#2 \arrowtypec#3}
\def\settoheight#1#2{\setbox\@tempboxa\hbox{#2}#1\ht\@tempboxa\relax}%
\def\settodepth#1#2{\setbox\@tempboxa\hbox{#2}#1\dp\@tempboxa\relax}%
\def\settokens`#1`#2`#3`#4`{%
     \def\tokena{#1}\def\tokenb{#2}\def\tokenc{#3}\def\tokend{#4}}
\def\setsqparms[#1`#2`#3`#4;#5`#6]{%
\arrowtypea #1
\arrowtypeb #2
\arrowtypec #3
\arrowtyped #4
\width #5
\height #6
}
\def\setpos(#1,#2){\xpos=#1 \ypos#2}

\def\settriparms[#1`#2`#3;#4]{\settripairparms[#1`#2`#3`1`1;#4]}%

\def\settripairparms[#1`#2`#3`#4`#5;#6]{%
\arrowtypea #1
\arrowtypeb #2
\arrowtypec #3
\arrowtyped #4
\arrowtypee #5
\width #6
\height #6
}

\def\resetparms{\settripairparms[1`1`1`1`1;500]\width 500}

\resetparms

\def\mvector(#1,#2)#3{
\put(0,0){\vector(#1,#2){#3}}%
\put(0,0){\vector(#1,#2){26}}%
}
\def\evector(#1,#2)#3{{
\arrowlength #3
\put(0,0){\vector(#1,#2){\arrowlength}}%
\advance \arrowlength by-30
\put(0,0){\vector(#1,#2){\arrowlength}}%
}}

\def\horsize#1#2{%
\settowidth{\tempdimen}{$#2$}%
#1=\tempdimen
\divide #1 by\unitlength
}

\def\vertsize#1#2{%
\settoheight{\tempdimen}{$#2$}%
#1=\tempdimen
\settodepth{\tempdimen}{$#2$}%
\advance #1 by\tempdimen
\divide #1 by\unitlength
}

\def\putvector(#1,#2)(#3,#4)#5#6{{%
\ifnum3<\arrowtype
\putdashvector(#1,#2)(#3,#4)#5\arrowtype
\else
\ifnum\arrowtype<-3
\putdashvector(#1,#2)(#3,#4)#5\arrowtype
\else
\xpos=#1
\ypos=#2
\run=#3
\rise=#4
\arrowlength=#5
\ifnum \arrowtype<0
    \ifnum \run=0
        \advance \ypos by-\arrowlength
    \else
        \tempcounta \arrowlength
        \multiply \tempcounta by\rise
        \divide \tempcounta by\run
        \ifnum\run>0
            \advance \xpos by\arrowlength
            \advance \ypos by\tempcounta
        \else
            \advance \xpos by-\arrowlength
            \advance \ypos by-\tempcounta
        \fi
    \fi
    \multiply \arrowtype by-1
    \multiply \rise by-1
    \multiply \run by-1
\fi
\ifcase \arrowtype
\or \put(\xpos,\ypos){\vector(\run,\rise){\arrowlength}}%
\or \put(\xpos,\ypos){\mvector(\run,\rise)\arrowlength}%
\or \put(\xpos,\ypos){\evector(\run,\rise){\arrowlength}}%
\fi\fi\fi
}}

\def\putsplitvector(#1,#2)#3#4{
\xpos #1
\ypos #2
\arrowtype #4
\halflength #3
\arrowlength #3
\gap 140
\advance \halflength by-\gap
\divide \halflength by2
\ifnum\arrowtype>0
   \ifcase \arrowtype
   \or \put(\xpos,\ypos){\line(0,-1){\halflength}}%
       \advance\ypos by-\halflength
       \advance\ypos by-\gap
       \put(\xpos,\ypos){\vector(0,-1){\halflength}}%
   \or \put(\xpos,\ypos){\line(0,-1)\halflength}%
       \put(\xpos,\ypos){\vector(0,-1)3}%
       \advance\ypos by-\halflength
       \advance\ypos by-\gap
       \put(\xpos,\ypos){\vector(0,-1){\halflength}}%
   \or \put(\xpos,\ypos){\line(0,-1)\halflength}%
       \advance\ypos by-\halflength
       \advance\ypos by-\gap
       \put(\xpos,\ypos){\evector(0,-1){\halflength}}%
   \fi
\else \arrowtype=-\arrowtype
   \ifcase\arrowtype
   \or \advance \ypos by-\arrowlength
       \put(\xpos,\ypos){\line(0,1){\halflength}}%
       \advance\ypos by\halflength
       \advance\ypos by\gap
       \put(\xpos,\ypos){\vector(0,1){\halflength}}%
   \or \advance \ypos by-\arrowlength
       \put(\xpos,\ypos){\line(0,1)\halflength}%
       \put(\xpos,\ypos){\vector(0,1)3}%
       \advance\ypos by\halflength
       \advance\ypos by\gap
       \put(\xpos,\ypos){\vector(0,1){\halflength}}%
   \or \advance \ypos by-\arrowlength
       \put(\xpos,\ypos){\line(0,1)\halflength}%
       \advance\ypos by\halflength
       \advance\ypos by\gap
       \put(\xpos,\ypos){\evector(0,1){\halflength}}%
   \fi
\fi
}

\def\putmorphism(#1)(#2,#3)[#4`#5`#6]#7#8#9{{%
\run #2
\rise #3
\ifnum\rise=0
  \puthmorphism(#1)[#4`#5`#6]{#7}{#8}#9%
\else\ifnum\run=0
  \putvmorphism(#1)[#4`#5`#6]{#7}{#8}#9%
\else
\setpos(#1)%
\arrowlength #7
\arrowtype #8
\ifnum\run=0
\else\ifnum\rise=0
\else
\ifnum\run>0
    \coefa=1
\else
   \coefa=-1
\fi
\ifnum\arrowtype>0
   \coefb=0
   \coefc=-1
\else
   \coefb=\coefa
   \coefc=1
   \arrowtype=-\arrowtype
\fi
\width=2
\multiply \width by\run
\divide \width by\rise
\ifnum \width<0  \width=-\width\fi
\advance\width by60
\if l#9 \width=-\width\fi
\putbox(\xpos,\ypos){#4}
{\multiply \coefa by\arrowlength
\advance\xpos by\coefa
\multiply \coefa by\rise
\divide \coefa by\run
\advance \ypos by\coefa
\putbox(\xpos,\ypos){#5} }%
{\multiply \coefa by\arrowlength
\divide \coefa by2
\advance \xpos by\coefa
\advance \xpos by\width
\multiply \coefa by\rise
\divide \coefa by\run
\advance \ypos by\coefa
\if l#9%
   \putrbox(\xpos,\ypos){#6}%
\else\if r#9%
   \putlbox(\xpos,\ypos){#6}%
\fi\fi }%
{\multiply \rise by-\coefc
\multiply \run by-\coefc
\multiply \coefb by\arrowlength
\advance \xpos by\coefb
\multiply \coefb by\rise
\divide \coefb by\run
\advance \ypos by\coefb
\multiply \coefc by70
\advance \ypos by\coefc
\multiply \coefc by\run
\divide \coefc by\rise
\advance \xpos by\coefc
\multiply \coefa by140
\multiply \coefa by\run
\divide \coefa by\rise
\advance \arrowlength by\coefa
\ifcase\arrowtype
\or \put(\xpos,\ypos){\vector(\run,\rise){\arrowlength}}%
\or \put(\xpos,\ypos){\mvector(\run,\rise){\arrowlength}}%
\or \put(\xpos,\ypos){\evector(\run,\rise){\arrowlength}}%
\fi}\fi\fi\fi\fi}}

\newcount\numbdashes \newcount\lengthdash \newcount\increment

\def\howmanydashes{
\numbdashes=\arrowlength \lengthdash=40
\divide\numbdashes by \lengthdash
\lengthdash=\arrowlength
\divide\lengthdash by \numbdashes
\increment=\lengthdash
\multiply\lengthdash by 3
\divide\lengthdash by 5
}

\def\putdashvector(#1)(#2,#3)#4#5{%
\ifnum#3=0 \putdashhvector(#1){#4}#5
\else
\ifnum#2=0
\putdashvvector(#1){#4}#5\fi\fi}

\def\putdashhvector(#1,#2)#3#4{{%
\arrowlength=#3 \howmanydashes
\multiput(#1,#2)(\increment,0){\numbdashes}%
{\vrule height .4pt width \lengthdash\unitlength}
\arrowtype=#4 \xpos=#1
\ifnum\arrowtype<0 \advance\arrowtype by 7 \fi
\ifcase\arrowtype
\or \advance\xpos by 10
    \put(\xpos,#2){\vector(-1,0){\lengthdash}}
    \advance\xpos by 40
    \put(\xpos,#2){\vector(-1,0){\lengthdash}}
\or \advance \xpos by 10
    \put(\xpos,#2){\vector(-1,0){\lengthdash}}
    \advance\xpos by  \arrowlength
    \advance\xpos by  -50
    \put(\xpos,#2){\vector(-1,0){\lengthdash}}
\or \advance\xpos by 10
    \put(\xpos,#2){\vector(-1,0){\lengthdash}}
\or \advance\xpos by \arrowlength
    \advance\xpos by -\lengthdash
    \put(\xpos,#2){\vector(1,0){\lengthdash}}
\or {\advance\xpos by 10
    \put(\xpos,#2){\vector(1,0){\lengthdash}}}
    \advance\xpos by \arrowlength
    \advance\xpos by -\lengthdash
    \put(\xpos,#2){\vector(1,0){\lengthdash}}
\or \advance\xpos by \arrowlength
    \advance\xpos by -\lengthdash
    \put(\xpos,#2){\vector(1,0){\lengthdash}}
    \advance\xpos by -40
    \put(\xpos,#2){\vector(1,0){\lengthdash}}
   \fi
}}

\def\putdashvvector(#1,#2)#3#4{{%
\arrowlength=#3 \howmanydashes
\ypos=#2 \advance\ypos by -\arrowlength
\multiput(#1,#2)(0,\increment){\numbdashes}%
    {\vrule width .4pt height \lengthdash\unitlength}
\arrowtype=#4 \ypos=#2
\ifnum\arrowtype<0 \advance\arrowtype by 7 \fi
\ifcase\arrowtype
\or \advance\ypos by \arrowlength \advance\ypos by -40
    \put(#1,\ypos){\vector(0,1){\lengthdash}}
    \advance\ypos by -40
    \put(#1,\ypos){\vector(0,1){\lengthdash}}
\or \advance\ypos by 10
    \put(#1,\ypos){\vector(0,1){\lengthdash}}
    \advance\ypos by \arrowlength \advance\ypos by -40
    \put(#1,\ypos){\vector(0,1){\lengthdash}}
\or \advance\ypos by \arrowlength \advance\ypos by -40
    \put(#1,\ypos){\vector(0,1){\lengthdash}}
\or \advance\ypos by 10
    \put(#1,\ypos){\vector(0,-1){\lengthdash}}
\or \advance\ypos by 10
    \put(#1,\ypos){\vector(0,-1){\lengthdash}}
    \advance\ypos by \arrowlength \advance\ypos by -40
    \put(#1,\ypos){\vector(0,-1){\lengthdash}}
\or \advance\ypos by 10
    \put(#1,\ypos){\vector(0,-1){\lengthdash}}
    \advance\ypos by 40
    \put(#1,\ypos){\vector(0,-1){\lengthdash}}
\fi
}}

\def\puthmorphism(#1,#2)[#3`#4`#5]#6#7#8{{%
\xpos #1
\ypos #2
\width #6
\arrowlength #6
\arrowtype=#7
\putbox(\xpos,\ypos){#3\vphantom{#4}}%
{\advance \xpos by\arrowlength
\putbox(\xpos,\ypos){\vphantom{#3}#4}}%
\horsize{\tempcounta}{#3}%
\horsize{\tempcountb}{#4}%
\divide \tempcounta by2
\divide \tempcountb by2
\advance \tempcounta by30
\advance \tempcountb by30
\advance \xpos by\tempcounta
\advance \arrowlength by-\tempcounta
\advance \arrowlength by-\tempcountb
\putvector(\xpos,\ypos)(1,0)\arrowlength\arrowtype
\divide \arrowlength by2
\advance \xpos by\arrowlength
\vertsize{\tempcounta}{#5}%
\divide\tempcounta by2
\advance \tempcounta by20
\if a#8 %
   \advance \ypos by\tempcounta
   \putbox(\xpos,\ypos){#5}%
\else
   \advance \ypos by-\tempcounta
   \putbox(\xpos,\ypos){#5}%
\fi}}

\def\putvmorphism(#1,#2)[#3`#4`#5]#6#7#8{{%
\xpos #1
\ypos #2
\arrowlength #6
\arrowtype #7
\settowidth{\xlen}{$#5$}%
\putbox(\xpos,\ypos){#3}%
{\advance \ypos by-\arrowlength
\putbox(\xpos,\ypos){#4}}%
{\advance\arrowlength by-140
\advance \ypos by-70
\ifdim\xlen>0pt
   \if m#8%
      \putsplitvector(\xpos,\ypos)\arrowlength\arrowtype
   \else
   \putvector(\xpos,\ypos)(0,-1)\arrowlength\arrowtype
   \fi
\else
   \putvector(\xpos,\ypos)(0,-1)\arrowlength\arrowtype
\fi}%
\ifdim\xlen>0pt
   \divide \arrowlength by2
   \advance\ypos by-\arrowlength
   \if l#8%
      \advance \xpos by-40
      \putrbox(\xpos,\ypos){#5}%
   \else\if r#8%
      \advance \xpos by40
      \putlbox(\xpos,\ypos){#5}%
   \else
      \putbox(\xpos,\ypos){#5}%
   \fi\fi
\fi
}}

\def\putsquarep<#1>(#2)[#3;#4`#5`#6`#7]{{%
\setsqparms[#1]%
\setpos(#2)%
\settokens`#3`%
\puthmorphism(\xpos,\ypos)[\tokenc`\tokend`{#7}]{\width}{\arrowtyped}b%
\advance\ypos by \height
\puthmorphism(\xpos,\ypos)[\tokena`\tokenb`{#4}]{\width}{\arrowtypea}a%
\putvmorphism(\xpos,\ypos)[``{#5}]{\height}{\arrowtypeb}l%
\advance\xpos by \width
\putvmorphism(\xpos,\ypos)[``{#6}]{\height}{\arrowtypec}r%
}}

\def\putsquare{\@ifnextchar <{\putsquarep}{\putsquarep%
   <\arrowtypea`\arrowtypeb`\arrowtypec`\arrowtyped;\width`\height>}}
\def\squaree{\@ifnextchar< {\squarep}{\squarep
   <\arrowtypea`\arrowtypeb`\arrowtypec`\arrowtyped;\width`\height>}}
\def\squarep<#1>[#2`#3`#4`#5;#6`#7`#8`#9]{{
\setsqparms[#1]
\diagram
\putsquarep<\arrowtypea`\arrowtypeb`\arrowtypec`
\arrowtyped;\width`\height>
(0,0)[#2`#3`#4`{#5};#6`#7`#8`{#9}]
\enddiagram
}}                                                 
\def\putptrianglep<#1>(#2,#3)[#4`#5`#6;#7`#8`#9]{{%
\settriparms[#1]%
\xpos=#2 \ypos=#3
\advance\ypos by \height
\puthmorphism(\xpos,\ypos)[#4`#5`{#7}]{\height}{\arrowtypea}a%
\putvmorphism(\xpos,\ypos)[`#6`{#8}]{\height}{\arrowtypeb}l%
\advance\xpos by\height
\putmorphism(\xpos,\ypos)(-1,-1)[``{#9}]{\height}{\arrowtypec}r%
}}

\def\putptriangle{\@ifnextchar <{\putptrianglep}{\putptrianglep
   <\arrowtypea`\arrowtypeb`\arrowtypec;\height>}}
\def\ptriangle{\@ifnextchar <{\ptrianglep}{\ptrianglep
   <\arrowtypea`\arrowtypeb`\arrowtypec;\height>}}
\def\ptrianglep<#1>[#2`#3`#4;#5`#6`#7]{{
\settriparms[#1]
\diagram
\putptrianglep<\arrowtypea`\arrowtypeb`
\arrowtypec;\height>
(0,0)[#2`#3`#4;#5`#6`{#7}]
\enddiagram
}}                                            

\def\putqtrianglep<#1>(#2,#3)[#4`#5`#6;#7`#8`#9]{{%
\settriparms[#1]%
\xpos=#2 \ypos=#3
\advance\ypos by\height
\puthmorphism(\xpos,\ypos)[#4`#5`{#7}]{\height}{\arrowtypea}a%
\putmorphism(\xpos,\ypos)(1,-1)[``{#8}]{\height}{\arrowtypeb}l%
\advance\xpos by\height
\putvmorphism(\xpos,\ypos)[`#6`{#9}]{\height}{\arrowtypec}r%
}}

\def\putqtriangle{\@ifnextchar <{\putqtrianglep}{\putqtrianglep
   <\arrowtypea`\arrowtypeb`\arrowtypec;\height>}}
\def\qtriangle{\@ifnextchar <{\qtrianglep}{\qtrianglep
   <\arrowtypea`\arrowtypeb`\arrowtypec;\height>}}
\def\qtrianglep<#1>[#2`#3`#4;#5`#6`#7]{{
\settriparms[#1]
\width=\height                                
\diagram
\putqtrianglep<\arrowtypea`\arrowtypeb`
\arrowtypec;\height>
(0,0)[#2`#3`#4;#5`#6`{#7}]
\enddiagram
}}

\def\putdtrianglep<#1>(#2,#3)[#4`#5`#6;#7`#8`#9]{{%
\settriparms[#1]%
\xpos=#2 \ypos=#3
\puthmorphism(\xpos,\ypos)[#5`#6`{#9}]{\height}{\arrowtypec}b%
\advance\xpos by \height \advance\ypos by\height
\putmorphism(\xpos,\ypos)(-1,-1)[``{#7}]{\height}{\arrowtypea}l%
\putvmorphism(\xpos,\ypos)[#4``{#8}]{\height}{\arrowtypeb}r%
}}

\def\putdtriangle{\@ifnextchar <{\putdtrianglep}{\putdtrianglep
   <\arrowtypea`\arrowtypeb`\arrowtypec;\height>}}
\def\dtriangle{\@ifnextchar <{\dtrianglep}{\dtrianglep
   <\arrowtypea`\arrowtypeb`\arrowtypec;\height>}}
\def\dtrianglep<#1>[#2`#3`#4;#5`#6`#7]{{
\settriparms[#1]
\width=\height                                
\diagram
\putdtrianglep<\arrowtypea`\arrowtypeb`
\arrowtypec;\height>
(0,0)[#2`#3`#4;#5`#6`{#7}]
\enddiagram
}}

\def\putbtrianglep<#1>(#2,#3)[#4`#5`#6;#7`#8`#9]{{%
\settriparms[#1]%
\xpos=#2 \ypos=#3
\puthmorphism(\xpos,\ypos)[#5`#6`{#9}]{\height}{\arrowtypec}b%
\advance\ypos by\height
\putmorphism(\xpos,\ypos)(1,-1)[``{#8}]{\height}{\arrowtypeb}r%
\putvmorphism(\xpos,\ypos)[#4``{#7}]{\height}{\arrowtypea}l%
}}

\def\putbtriangle{\@ifnextchar <{\putbtrianglep}{\putbtrianglep
   <\arrowtypea`\arrowtypeb`\arrowtypec;\height>}}
\def\btriangle{\@ifnextchar <{\btrianglep}{\btrianglep
   <\arrowtypea`\arrowtypeb`\arrowtypec;\height>}}
\def\btrianglep<#1>[#2`#3`#4;#5`#6`#7]{{
\settriparms[#1]
\width=\height                               
\diagram
\putbtrianglep<\arrowtypea`\arrowtypeb`
\arrowtypec;\height>
(0,0)[#2`#3`#4;#5`#6`{#7}]
\enddiagram
}}

\def\putAtrianglep<#1>(#2,#3)[#4`#5`#6;#7`#8`#9]{{%
\settriparms[#1]%
\xpos=#2 \ypos=#3
{\multiply \height by2
\puthmorphism(\xpos,\ypos)[#5`#6`{#9}]{\height}{\arrowtypec}b}%
\advance\xpos by\height \advance\ypos by\height
\putmorphism(\xpos,\ypos)(-1,-1)[#4``{#7}]{\height}{\arrowtypea}l%
\putmorphism(\xpos,\ypos)(1,-1)[``{#8}]{\height}{\arrowtypeb}r%
}}

\def\putAtriangle{\@ifnextchar <{\putAtrianglep}{\putAtrianglep
   <\arrowtypea`\arrowtypeb`\arrowtypec;\height>}}
\def\Atriangle{\@ifnextchar <{\Atrianglep}{\Atrianglep
   <\arrowtypea`\arrowtypeb`\arrowtypec;\height>}}
\def\Atrianglep<#1>[#2`#3`#4;#5`#6`#7]{{
\settriparms[#1]
\width=\height                                     
\diagram
\putAtrianglep<\arrowtypea`\arrowtypeb`
\arrowtypec;\height>
(0,0)[#2`#3`#4;#5`#6`{#7}]
\enddiagram
}}

\def\putAtrianglepairp<#1>(#2)[#3;#4`#5`#6`#7`#8]{{%
\settripairparms[#1]%
\setpos(#2)%
\settokens`#3`%
\puthmorphism(\xpos,\ypos)[\tokenb`\tokenc`{#7}]{\height}{\arrowtyped}b%
\advance\xpos by\height
\puthmorphism(\xpos,\ypos)[\phantom{\tokenc}`\tokend`{#8}]%
{\height}{\arrowtypee}b%
\advance\ypos by\height
\putmorphism(\xpos,\ypos)(-1,-1)[\tokena``{#4}]{\height}{\arrowtypea}l%
\putvmorphism(\xpos,\ypos)[``{#5}]{\height}{\arrowtypeb}m%
\putmorphism(\xpos,\ypos)(1,-1)[``{#6}]{\height}{\arrowtypec}r%
}}

\def\putAtrianglepair{\@ifnextchar <{\putAtrianglepairp}{\putAtrianglepairp%
   <\arrowtypea`\arrowtypeb`\arrowtypec`\arrowtyped`\arrowtypee;\height>}}
\def\Atrianglepair{\@ifnextchar <{\Atrianglepairp}{\Atrianglepairp%
   <\arrowtypea`\arrowtypeb`\arrowtypec`\arrowtyped`\arrowtypee;\height>}}

\def\Atrianglepairp<#1>[#2;#3`#4`#5`#6`#7]{{
\settripairparms[#1]
\settokens`#2`
\width=\height                                
\diagram
\putAtrianglepairp                            
<\arrowtypea`\arrowtypeb`\arrowtypec`
\arrowtyped`\arrowtypee;\height>
(0,0)[{#2};#3`#4`#5`#6`{#7}]
\enddiagram
}}

\def\putVtrianglep<#1>(#2,#3)[#4`#5`#6;#7`#8`#9]{{%
\settriparms[#1]%
\xpos=#2 \ypos=#3
\advance\ypos by\height
{\multiply\height by2
\puthmorphism(\xpos,\ypos)[#4`#5`{#7}]{\height}{\arrowtypea}a}%
\putmorphism(\xpos,\ypos)(1,-1)[`#6`{#8}]{\height}{\arrowtypeb}l%
\advance\xpos by\height
\advance\xpos by\height
\putmorphism(\xpos,\ypos)(-1,-1)[``{#9}]{\height}{\arrowtypec}r%
}}

\def\putVtriangle{\@ifnextchar <{\putVtrianglep}{\putVtrianglep
   <\arrowtypea`\arrowtypeb`\arrowtypec;\height>}}
\def\Vtriangle{\@ifnextchar <{\Vtrianglep}{\Vtrianglep
   <\arrowtypea`\arrowtypeb`\arrowtypec;\height>}}
\def\Vtrianglep<#1>[#2`#3`#4;#5`#6`#7]{{
\settriparms[#1]
\width=\height                                 
\diagram
\putVtrianglep<\arrowtypea`\arrowtypeb`
\arrowtypec;\height>
(0,0)[#2`#3`#4;#5`#6`{#7}]
\enddiagram
}}

\def\putVtrianglepairp<#1>(#2)[#3;#4`#5`#6`#7`#8]{{
\settripairparms[#1]%
\setpos(#2)%
\settokens`#3`%
\advance\ypos by\height
\putmorphism(\xpos,\ypos)(1,-1)[`\tokend`{#6}]{\height}{\arrowtypec}l%
\puthmorphism(\xpos,\ypos)[\tokena`\tokenb`{#4}]{\height}{\arrowtypea}a%
\advance\xpos by\height
\puthmorphism(\xpos,\ypos)[\phantom{\tokenb}`\tokenc`{#5}]%
{\height}{\arrowtypeb}a%
\putvmorphism(\xpos,\ypos)[``{#7}]{\height}{\arrowtyped}m%
\advance\xpos by\height
\putmorphism(\xpos,\ypos)(-1,-1)[``{#8}]{\height}{\arrowtypee}r%
}}

\def\putVtrianglepair{\@ifnextchar <{\putVtrianglepairp}{\putVtrianglepairp%
    <\arrowtypea`\arrowtypeb`\arrowtypec`\arrowtyped`\arrowtypee;\height>}}
\def\Vtrianglepair{\@ifnextchar <{\Vtrianglepairp}{\Vtrianglepairp%
    <\arrowtypea`\arrowtypeb`\arrowtypec`\arrowtyped`\arrowtypee;\height>}}
\def\Vtrianglepairp<#1>[#2;#3`#4`#5`#6`#7]{{
\settripairparms[#1]
\settokens`#2`
\diagram
\putVtrianglepairp                             
<\arrowtypea`\arrowtypeb`\arrowtypec`
\arrowtyped`\arrowtypee;\height>
(0,0)[{#2};#3`#4`#5`#6`{#7}]
\enddiagram
}}

\def\putCtrianglep<#1>(#2,#3)[#4`#5`#6;#7`#8`#9]{{%
\settriparms[#1]%
\xpos=#2 \ypos=#3
\advance\ypos by\height
\putmorphism(\xpos,\ypos)(1,-1)[``{#9}]{\height}{\arrowtypec}l%
\advance\xpos by\height
\advance\ypos by\height
\putmorphism(\xpos,\ypos)(-1,-1)[#4`#5`{#7}]{\height}{\arrowtypea}l%
{\multiply\height by 2
\putvmorphism(\xpos,\ypos)[`#6`{#8}]{\height}{\arrowtypeb}r}%
}}

\def\putCtriangle{\@ifnextchar <{\putCtrianglep}{\putCtrianglep
    <\arrowtypea`\arrowtypeb`\arrowtypec;\height>}}
\def\Ctriangle{\@ifnextchar <{\Ctrianglep}{\Ctrianglep
    <\arrowtypea`\arrowtypeb`\arrowtypec;\height>}}
\def\Ctrianglep<#1>[#2`#3`#4;#5`#6`#7]{{
\settriparms[#1]
\width=\height                               
\diagram
\putCtrianglep<\arrowtypea`\arrowtypeb`
\arrowtypec;\height>
(0,0)[#2`#3`#4;#5`#6`{#7}]
\enddiagram
}}                                           
\def\putDtrianglep<#1>(#2,#3)[#4`#5`#6;#7`#8`#9]{{%
\settriparms[#1]%
\xpos=#2 \ypos=#3
\advance\xpos by\height \advance\ypos by\height
\putmorphism(\xpos,\ypos)(-1,-1)[``{#9}]{\height}{\arrowtypec}r%
\advance\xpos by-\height \advance\ypos by\height
\putmorphism(\xpos,\ypos)(1,-1)[`#5`{#8}]{\height}{\arrowtypeb}r%
{\multiply\height by 2
\putvmorphism(\xpos,\ypos)[#4`#6`{#7}]{\height}{\arrowtypea}l}%
}}

\def\putDtriangle{\@ifnextchar <{\putDtrianglep}{\putDtrianglep
    <\arrowtypea`\arrowtypeb`\arrowtypec;\height>}}
\def\Dtriangle{\@ifnextchar <{\Dtrianglep}{\Dtrianglep
   <\arrowtypea`\arrowtypeb`\arrowtypec;\height>}}
\def\Dtrianglep<#1>[#2`#3`#4;#5`#6`#7]{{
\settriparms[#1]
\width=\height                              
\diagram
\putDtrianglep<\arrowtypea`\arrowtypeb`
\arrowtypec;\height>
(0,0)[#2`#3`#4;#5`#6`{#7}]
\enddiagram
}}                                          
\def\setrecparms[#1`#2]{\width=#1 \height=#2}%

\def\recursep<#1`#2>[#3;#4`#5`#6`#7`#8]{{\m@th
\width=#1 \height=#2
\settokens`#3`
\settowidth{\tempdimen}{$\tokena$}
\ifdim\tempdimen=0pt
  \savebox{\tempboxa}{\hbox{$\tokenb$}}%
  \savebox{\tempboxb}{\hbox{$\tokend$}}%
  \savebox{\tempboxc}{\hbox{$#6$}}%
\else
  \savebox{\tempboxa}{\hbox{$\hbox{$\tokena$}\times\hbox{$\tokenb$}$}}%
  \savebox{\tempboxb}{\hbox{$\hbox{$\tokena$}\times\hbox{$\tokend$}$}}%
  \savebox{\tempboxc}{\hbox{$\hbox{$\tokena$}\times\hbox{$#6$}$}}%
\fi
\ypos=\height
\divide\ypos by 2
\xpos=\ypos
\advance\xpos by \width
\bfig
\putCtrianglep<-1`1`1;\ypos>(0,0)[`\tokenc`;#5`#6`{#7}]%
\puthmorphism(\ypos,0)[\tokend`\usebox{\tempboxb}`{#8}]{\width}{-1}b%
\puthmorphism(\ypos,\height)[\tokenb`\usebox{\tempboxa}`{#4}]{\width}{-1}a%
\advance\ypos by \width
\putvmorphism(\ypos,\height)[``\usebox{\tempboxc}]{\height}1r%
\efig
}}

\def\recurse{\@ifnextchar <{\recursep}{\recursep<\width`\height>}}

\def\puttwohmorphisms(#1,#2)[#3`#4;#5`#6]#7#8#9{{%
\puthmorphism(#1,#2)[#3`#4`]{#7}0a
\ypos=#2
\advance\ypos by 20
\puthmorphism(#1,\ypos)[\phantom{#3}`\phantom{#4}`#5]{#7}{#8}a
\advance\ypos by -40
\puthmorphism(#1,\ypos)[\phantom{#3}`\phantom{#4}`#6]{#7}{#9}b
}}

\def\puttwovmorphisms(#1,#2)[#3`#4;#5`#6]#7#8#9{{%
\putvmorphism(#1,#2)[#3`#4`]{#7}0a
\xpos=#1
\advance\xpos by -20
\putvmorphism(\xpos,#2)[\phantom{#3}`\phantom{#4}`#5]{#7}{#8}l
\advance\xpos by 40
\putvmorphism(\xpos,#2)[\phantom{#3}`\phantom{#4}`#6]{#7}{#9}r
}}

\def\puthcoequalizer(#1)[#2`#3`#4;#5`#6`#7]#8#9{{%
\setpos(#1)%
\puttwohmorphisms(\xpos,\ypos)[#2`#3;#5`#6]{#8}11%
\advance\xpos by #8
\puthmorphism(\xpos,\ypos)[\phantom{#3}`#4`#7]{#8}1{#9}
}}

\def\putvcoequalizer(#1)[#2`#3`#4;#5`#6`#7]#8#9{{%
\setpos(#1)%
\puttwovmorphisms(\xpos,\ypos)[#2`#3;#5`#6]{#8}11%
\advance\ypos by -#8
\putvmorphism(\xpos,\ypos)[\phantom{#3}`#4`#7]{#8}1{#9}
}}

\def\putthreehmorphisms(#1)[#2`#3;#4`#5`#6]#7(#8)#9{{%
\setpos(#1) \settypes(#8)
\if a#9 %
     \vertsize{\tempcounta}{#5}%
     \vertsize{\tempcountb}{#6}%
     \ifnum \tempcounta<\tempcountb \tempcounta=\tempcountb \fi
\else
     \vertsize{\tempcounta}{#4}%
     \vertsize{\tempcountb}{#5}%
     \ifnum \tempcounta<\tempcountb \tempcounta=\tempcountb \fi
\fi
\advance \tempcounta by 60
\puthmorphism(\xpos,\ypos)[#2`#3`#5]{#7}{\arrowtypeb}{#9}
\advance\ypos by \tempcounta
\puthmorphism(\xpos,\ypos)[\phantom{#2}`\phantom{#3}`#4]{#7}{\arrowtypea}{#9}
\advance\ypos by -\tempcounta \advance\ypos by -\tempcounta
\puthmorphism(\xpos,\ypos)[\phantom{#2}`\phantom{#3}`#6]{#7}{\arrowtypec}{#9}
}}

\def\setarrowtoks[#1`#2`#3`#4`#5`#6]{%
\def\toka{#1}
\def\tokb{#2}
\def\tokc{#3}
\def\tokd{#4}
\def\toke{#5}
\def\tokf{#6}
}
\def\hex{\@ifnextchar <{\hexp}{\hexp<1000`400>}}
\def\hexp<#1`#2>[#3`#4`#5`#6`#7`#8;#9]{%
\setarrowtoks[#9]
\yext=#2 \advance \yext by #2
\xext=#1 \advance\xext by \yext
\bfig
\putCtriangle<-1`0`1;#2>(0,0)[`#5`;\tokb``\tokd]
\xext=#1 \yext=#2 \advance \yext by #2
\putsquare<1`0`0`1;\xext`\yext>(#2,0)[#3`#4`#7`#8;\toka```\tokf]
\advance \xext by #2
\putDtriangle<0`1`-1;#2>(\xext,0)[`#6`;`\tokc`\toke]
\efig
}
\makeatother

\usepackage{epic,eepic}

\def\pSkip{\vskip 1.5mm \noindent}

\numberwithin{equation}{section}

\begin{document}


\title[Tropical Arithmetic and  Tropical Matrix Algebra]
{Tropical Arithmetic and  Tropical Matrix Algebra}

\author{Zur Izhakian}\thanks{The author has been supported by the
Chateaubriand scientific post-doctorate fellowships, Ministry of
Science, French Government, 2007-2008.}
\address{Department of Mathematics, Bar-Ilan University, Ramat-Gan 52900,
Israel} \address{ \vskip -6mm CNRS et Universit´e Denis Diderot
(Paris 7), 175, rue du Chevaleret 75013 Paris, France}
\email{zzur@math.biu.ac.il, zzur@post.tau.ac.il}


\date{Febuary 2008}




\begin{abstract}
This paper introduces a new structure of commutative semiring,
generalizing the tropical semiring, and having an arithmetic that
modifies the standard tropical operations, i.e. summation and
maximum. Although our framework is combinatorial, notions of
regularity and invertibility arise naturally for matrices over
this semiring; we show that a tropical matrix is invertible if and
only if it is regular.
\end{abstract}

\maketitle



\section*{Introduction}\label{sec:Introduction}
Traditionally,  researchers have been able to  frame mathematical
theories using formal structures provided by algebra; geometry is
often a source for interesting phenomena in the core of these
theories. The semiring structure introduced in this paper emerges
from the combinatorics within max-plus algebra and its
corresponding polyhedral geometry, called tropical geometry.
Although our ground structure is a semiring, much of the theory of
standard commutative algebra can be formulated on this semiring,
leading to application in combinatorics, semigroup theory,
polynomials algebra, and algebraic geometry.

Tropical mathematics takes place over the \bfem{tropical semiring}
$(\Real \cup \{ \tUniS \},\max, +)$, the real numbers equipped
with the operations of maximum and summation, respectively,
addition and multiplication \cite{Gathmann:0601322,IMS,RST}, and
it interacts with a number of fields of study including algebraic
geometry, polyhedral geometry, commutative algebra, and
combinatorics. Polyhedral complexes, resembling algebraic
varieties over a field with real non-archimedean valuation, are
the main objects of the tropical geometry, where their geometric
combinatorial structure is a maximal degeneration of a complex
structure on a manifold.

Over the past few years, much effort has been invested in the
attempt to characterize a tropical analogous to classical linear
algebra,  \cite{Develin2003,Kim:130853,RST}, and to determine
connections between the classical and the tropical worlds
\cite{M4,Shustin2005,Speyer4218}. Despite the progress that has
been achieved in these tropical studies, some fundamental issues
have not been settled yet; the idempotency of addition in $(\Real
\cup \{ \tUniS \},\max, +)$ is maybe one of the main reasons for
that. Addressing this reason, and other algebro-geometric needs,
our goals are:
\renewcommand{\labelenumi}{(\alph{enumi})}

\begin{enumerate}
\item  Introducing a new structure of a partial
idempotent semiring having its own arithmetic that generalizes the
max-plus arithmetic and also carries a tropical geometric meaning;
\pSkip
      \item  Presenting a novel approach for a theory of matrix algebra over partial idempotent semirings that includes
 notions of regularity and semigroup invertibility, analogous as possible to that of matrices over fields.
  \end{enumerate}
  \renewcommand{\labelenumi}{(\arabic{enumi})}
The latter goal is central issue in the study of Green's relations
over semigroups and is essential toward developing a linear
representations of semigroups. Our new approach answers these
goals and paves a way to treat other needs like having a notions
of linear dependency and rank.

Our new structure, which we call \bfem{extended tropical
semiring}, is built on the disjoint union of two copies of
$\Real$, denoted $\Real$ and $\tUnit$, together with the formal
element $\tUniS$ that serves as the gluing point of $\Real$ and
$\tUnit$. Thus,
$$\Trop := \Real \ \cup \ \{ \tUniS \} \  \cup \ \tUnit$$ is
provided with an order, $\prec$, extending the usual order on
$\Real$, and endowed with the addition $\TrS$ and the
multiplication $\TrP$ that modify the familiar operations $\max$
and $+$. By this setting, $\tSR$ has the structure of a
commutative semiring, $\TrS$ is idempotent only on $\tUnit$, and
$\tSR$ allows to define a homomorphic relation to a field with
real non-archimedean valuation. From the point of view of
algebraic geometry, $\TrS$ encodes an additive multiplicity that
enables to define tropical algebraic sets in a natural manner.

The second part of the paper focuses mainly on introducing a
theory of matrix algebra over $\tSR$, reassembling the classical
theory of matrices over fields, that includes notions of
regularity and invertibility in a natural way with the following
relation:
\rTheorem{\ref{thm:inverseMatrix}}{A tropical matrix is pseudo
invertible if and only if it is tropically regular.}
\noindent We provide also an explicit characterization of the
pseudo inverse matrix $\itA$ of a regular matrix  $\tA$, which
turns out to be similar to that of the classical theory.
Concerning semigroup theory, we show that the monoid $M_n(\Trop)$
of matrices over $\tSR$ can be related to as an E-dense monoid in
which our invertibility suits E-denseness, that is the products
$\tA \itA$ and $\itA \tA$ are idempotent matrices
\cite{margolisEdense}.

\pSkip
\textbf{\emph{Acknowledgement}}: The author would like to thank
\emph{Prof. Eugenii Shustin} for his invaluable help. I'm deeply
grateful him for his support and the fertile discussions we had.

A part of this work was done during the author's stay at the
\emph{Max-Planck-Institut f\"ur Mathematik (Bonn)}. The author is
very grateful to MPI for the hospitality and excellent work
conditions.


\section{Extended Tropical Arithmetic -- A New
Approach}\label{chap:TropicalArithmetic}


With two goals in minds, geometrically  and algebraically derived,
our objective is to introduce a new concept of idempotent semiring
extensions, applied here to the classical tropical semiring
$(\Real \cup \{ \tUniS \}, \max,+)$, including also the relation
to non-Archimedean fields with real valuations. Although related
topics have been discussed earlier for $(\Real \cup \{ \tUniS \},
\max,+)$, cf.
\cite{Butkovic2003,Cuninghame2004,Develin2003,SpeyerSturmfels2004},
in this paper we use a different approach implemented on a
semiring structure having a modified arithmetic. We open by
describing the standard tropical framework, then we present the
basics of our new concept and the associated semiring structure.

\subsection{The tropical semiring}\label{sec:Max-PlusArithmetics}
 Tropical mathematics
is the mathematics over idempotent semirings,  the \bfem{tropical
semiring} is usually taken to be $(\Real \cup \{ \tUniS\}, \max, +
\;)$;  the real numbers together with the formal element $\tUniS$,
and with the operations of tropical addition and tropical
multiplication
$$a + b \ := \ \max\{a,b\} \ ,  \qquad \quad a \cdot b \ := \ a+b \ ,$$
cf. \cite{M4,RST}. We write $\eReal$ for $\Real\cup\{-\infty\}$
and equip $\eReal^* :=\Real$ with the Euclidean topology, assuming
that $\eReal$ is homeomorphic to $[0,\infty)$.
 The tropical semiring contains the \bfem{max-plus algebra}
\cite{Cuninghame2004,RST} and it emerges as a target of
non-Archimedean fields with real valuation; it is an idempotent
semiring, i.e. $a +a =a$, with the unit $\one_{\eReal} := 0$, and
the zero element $\zero_{\eReal} := \tUniS$.

Elements of the semiring $\eReal[\lm_1,\dots,\lm_n]$ are called
tropical polynomials in $n$ variables over $\eReal$ and are of the
form
\begin{equation}\label{eq:polyToFunc1} \tpF = \max_{\bfi \in\Omega}\{ \langle
\Lm ,\bfi \rangle+ \al_\bfi \} \ ,
\end{equation} where
 $\langle \, \cdot , \cdot \, \rangle$ stands for the
standard scalar product, $\Omega\subset\Int^n$ is a finite
nonempty set of points $\bfi =(i_1,\dots,i_n)$ with nonnegative
coordinates, $\al_\bfi \in\Real$ for all $\bfi \in\Omega$,  and
$\Lm = (\lm_1, \dots, \lm_n)$. The addition and multiplication of
polynomials  are defined according to the familiar law.

Any tropical polynomial $f \in
\eReal[\lm_1,\dots,\lm_n]\backslash\{-\infty\}$ determines a
piecewise linear convex function $\tfF:\Real^{(n)} \To \Real$.
But, in the tropical case, the map $\tpF\mapsto\tfF$ is not
injective, and one can reduce the polynomial semiring so as to
have only those elements needed to describe functions.

A tropical hypersurface is defined to be the domain of
non-differentiability, also called the corner locus, of
$\tilde{f}$ for some $f \in
\eReal[\lm_1,\dots,\lm_n]\backslash\{-\infty\}$. Therefore, points
of a tropical hypersurface can be specified as the points on which
the value of $\tilde{f}$ is attained by at least two monomials of
$\tpF$. This property is crucial for understanding the purpose of
incorporating additive multiplicities, it will be used later to
distinguish the corner locus from the other points of a domain.

One of our goals is to establish  a semiring structure that allows
one to realize (algebraically) the points of a corner locus as a
``zero'' locus of a polynomial; namely, to have the ability to
form algebraic sets. Therefore, we would like to have a structure
that not only provides the operation of maximum, but also encodes
an indication about its additive multiplicity. In other word, in
some sense, to ``resolve" the idempotency of $(\eReal, \max, + \
)$.

\begin{remark}\label{rem:nonAssociative} Indeed, to address  this goal,
one may suggest an alternative arithmetic that defines the
addition of two equal elements to be $\tUniS$, which we write as
$``a + a " = \tUniS$, and the addition of different elements to be
their maximum. We denote this structure as $(\eReal,``\max",+)$.
Unfortunately, this type of addition is not associative; for
example, for $b < a$ we have  $``\aab + \OP \aaa + \aaa \CP" =
``\aab + \OP \tUniS \CP" = \aab$ while $ ``\OP \aab + \aaa \CP +
\aaa" = ``(\aaa) + \aaa" = \tUniS$.
\end{remark}

 Our next development addresses this algebro-geometric issue;
 later we show that it also servers a solid base for developing a theory of matrix
 algebra over semirings
 that have the notions of  regularity and
 invertibility.
 \subsection{The extended tropical
semiring}\label{sec:TheExtendedTropicalArithmetic}

Roughly  speaking, the central idea of our new approach is a
generalization of $\eMaxPlus$ to a semiring structure having  a
partial idempotent addition that distinguishes between sums of
similar elements and sums of different elements. Set
theoretically, our semiring is composed from the disjoint union of
two copies of $\Real,$ denoted $\Real$ and $\tUnit$, which glued
along the formal element $-\infty$ to create the set $$  \Trop :=
\Real \ \cup \ \{ \tUniS\}  \ \cup  \ \tUnit \ .
$$ In what follows we denote
 the unions  $\Real \cup \{ \tUniS \}$ and $\tUnit
\cup \{ \tUniS \}$  respectively by $\eReal$ and $\etUnit$, write
$\Trop^{\times}$ for $\Trop \setminus \{ \tUniS \}$, and call the
elements of $\Real$ \bfem{reals}.

We use the generic notation that $a,b\in \Real$ for reals,
$\uuu{a},\uuu{b} \in \tUnit$ where  $a,b \in \Real$, and $x,y \in
\Trop$. Thus, $\Trop$ is provided with the following order $\prec$
extending the usual order on $\Real$:
\begin{axiom}\label{ax:tropRelations} The order $\prec$  on
$\Trop$ is defined as:
\begin{enumerate}
\item $-\infty \prec x,$ $\forall x \in \Trop^{\times}$; \pSkip

\item for any real numbers $a < b$, we have $a\prec
b,$ $a\prec b^\nu$, $a^\nu \prec b$, and $a^\nu \prec b^\nu;$
\pSkip

\item $a \prec a^\nu$ for all $a\in \Real.$ \pSkip
\end{enumerate}
\end{axiom}
\noindent One can verify that the corresponding partial order,
$\preceq$, holds only in the cases where both elements are in
$\Real$ or both are in $\tUnit$.

\begin{example} Assume $a< b < c $ are reals;  then
$$\tUniS  \ \prec \ a \ \prec \ \uuu{a} \ \prec \ b \ \prec \ \uuu{b} \ \prec \ c \ \prec \ \uuu{c}\  . $$

\end{example}

According to the rules of $\prec$, cf. Axiom
\ref{ax:tropRelations}, $\Trop$ is then endowed with the two
operations $\TrS$ and $\TrP$, addition and multiplication
respectively,  defined as below. (We use the notation
$\max_{\prec}$ to denote the maximum with respect to the order
$\prec$ \ .)

\begin{axiom}\label{ax:tropOperations} The laws of the extended
tropical arithmetic are:
\begin{enumerate}
\item $-\infty \oplus x = x \oplus -\infty = x$ for each $x \in
\Trop $; \pSkip

\item $x \TrS y =
\max_{\prec} \{x,y\}$ unless $x = y;$ \pSkip

\item $a\TrS a = a^\nu \TrS a^\nu =
a^\nu;$ \pSkip

\item $\tUniS \TrP x = x \TrP \tUniS = \tUniS$ for each $x \in \Trop
$; \pSkip

\item $a \TrP b
= a+b$ for all $a,b \in \Real;$ \pSkip

\item $ a^\nu \TrP b = a \TrP b ^\nu =
a^\nu \TrP b^\nu = (a+b)^\nu$ .\pSkip
\end{enumerate}
\end{axiom}
\noindent We call the triple $\tSR$ the \bfem{extended tropical
semiring}; later we show that $\tSR$ indeed have the structure of
commutative semiring with unit $\one_{\Trop} := 0$ and
$\zero_{\Trop} := \tUniS$.

 Recall that two preliminary essential demands have been
required on $\TrS$ , validity of associativity and,
simultaneously, differentiation between addition of similar reals
and addition of different reals. The first requirement is
satisfied by Axiom \ref{ax:tropRelations} (3) and Axiom
\ref{ax:tropOperations} (2); that is, for reals, we have the
following:
\begin{equation}\label{eq:exampleAcc}
   \aab \TrS \uuu{\aaa} = \aab \TrS \OP \aaa \TrS \aaa \CP \; \;
   \underset{\uparrow}{\overset{?}{=}} \; \;
   \OP \aab \TrS \aaa \CP \TrS \aaa =
\left\{ \begin{array}{lll}
  (\aaa) \TrS \aaa  & = \uuu{\aaa},& \aaa \succ \aab,  \\ [2mm]
  (\uuu{\aaa}) \TrS \aaa  & = \uuu{\aaa},& \aaa = \aab,  \\ [2mm]
  (\aab) \TrS \aaa  & = \aab, & \aab \succ \aaa,
\end{array}\right.
\end{equation}
 the equality is then derived from Axiom \ref{ax:tropOperations}
(2).

\begin{remark} $ $
\begin{enumerate}
    \item The addition $\TrS$ (in comparison  to that of $\eMaxPlus$) is not
idempotent, since $a \TrS a = \uuu{a}$; this is one of the main
aspects  of our approach. \pSkip
    \item $\etUnit$ is an ideal of $\Trop$ where  sometimes we want to
think about as a set of pseudo zeros, namely  consisting  of those
elements to be ignored. On the other hand, by Axiom
\ref{ax:tropOperations} (3), $\etUnit$ can be also realized as a
``shadow'' copy of $\Real$ whose elements carry additive
multiplicities $ > 1$, received as tropical sums of identical
reals. This view is important for understanding the linkage
between our arithmetic and the notion of tropicalization.

\end{enumerate}
\end{remark}

In the context of semigroups, both $(\Trop, \TrS)$ and
$(\Trop,\TrP)$ are monoids but not groups and thus, invertibility
is invalid for both $\TrS$ and $\TrP$ . Yet, for $\TrP$ , one can
talk about partial invertibility which is well defined on reals
only.
\begin{definition} The division, denoted $\iTrP$, of $x,y \in \Trop$, with $y \neq \tUniS$, is defined as
$x \iTrP y = x \TrP (-y)$, where $-y = \uuu{(-a)}$ when $y =
\uuu{a}$ .
\end{definition}
 Note that $\iTrP$  is not well defined over all $\Trop$,
 but suits our purpose. The
cancellation law, $ x \TrP y = x \TrP z \Rightarrow y = z$, does
not always hold; for example, the equality $\uuu{\aaa} \TrP \aab =
\uuu{\aaa} \TrP \uuu{\aab}$ does not satisfy  cancellation.

\begin{remark} The structure of $(\Trop, \TrS, \TrP)$ has been
formulated on two disjoint copies of $\Real$ with the modification
of the operations $\max$ and $+$ ; the same construction can be
performed for any idempotent semiring with the property $a + b \in
\{a,b \}$ and in particular for $(\Int, \max, + \ )$.
\end{remark}

\subsection{Properties of the extended tropical
arithmetic}\label{sec:PropertiesofTropicalOperators}
Having formulated extended tropical arithmetic, we address the its
basic properties. We describe only the main cases in detail;
therefore, the trivial cases involving $\tUniS$ are omitted. To
clarify the  exposition, sometimes, we treat the elements of
$\Real$ and $\tUnit$ separately.

\parSpc
\noindent \bfem{Commutativity:} Axiomatic (cf.~Axiom
    \ref{ax:tropOperations}).
\parSpc
\noindent    \bfem{Associativity:} By definition $\uuu{(a \TrS b)}
= \uuu{a} \TrS \uuu{b}$ and  $\uuu{(a \TrP b)} = \uuu{a} \TrP
    \uuu{b}$. Thus, for different elements in $\Trop$, the associativity of
$\TrS$ and $\TrP$ is clear by the associativity
    of  $\max$ and $+$ which also provides the associativity of
    $\TrP$ for all  $\Trop$. The case in which identical reals
    are involved has already been examined in \Ref{eq:exampleAcc}.
    For the case of two similar elements in $\tUnit$ we have:
   $$
   \begin{array}{lll}
     \OP \aaa \TrS \uuu{\aab} \CP \TrS \uuu{\aac} = &
      \left\{ \begin{array}{lll}
      \uuu{\aab} \TrS \uuu{\aac} = &  \uuu{(\aab \TrS \aac)}, & \aab \succeq \aaa, \\ [2mm]
      \aaa  \TrS \uuu{\aac} \searrow   &  & \aaa \succ \aab,
    \end{array}\right. \\ [4mm]
      &
       \hskip 2.3cm = \left\{ \begin{array}{lll}
      \uuu{ \aac},  &  & \aac \succeq \aaa,\aab, \\ [2mm]
      \aaa ,   &   & \aaa \succ \aab,\aac,
    \end{array}\right. \\ [6mm]
   \end{array}
   $$
   and
 $$ \begin{array}{lll}

      \aaa \TrS \OP \uuu{\aab} \TrS \uuu{\aac} \CP = &
      \aaa \TrS \uuu{( \aab  \TrS \aac)} =
       \left\{ \begin{array}{lll}
      \uuu{\aac}, &   & \aac \succeq \aaa, \aab,\\ [2mm]
      \uuu{\aab}, &   & \aab \succeq \aaa, \aac,\\ [2mm]
      \aaa,  &  & \aaa \succ \aab,\aac,
    \end{array}\right. \\
   \end{array}
   $$
   which have equal evaluations.
    (The other cases of compound expressions are obtained
     by the same way.)
\\
\parSpc \noindent
\bfem{Distributivity}: To verify distributivity of  $\TrP$ over
$\TrS$, for the  case when all elements are reals, write
$$
\begin{array}{lll}
      \aaa \TrP (\aab \TrS \aac)& =  & \left\{ \begin{array}{lll}
      \aaa \TrP \aab , &   & \aab \succ \aac,\\ [2mm]
      \aaa \TrP \uuu{\aab}, &   & \aab = \aac,\\ [2mm]
      \aaa \TrP \aac,  &  & \aac \succ \aab,
    \end{array}\right.
    \end{array}
    $$
and
$$    \begin{array}{lll}
  ( \aaa \TrP \aab) \TrS (\aaa \TrP \aac) & = &  \left\{ \begin{array}{lll}
      \aaa \TrP \aab , &   & \aab \succ \aac,\\ [2mm]
       \uuu{(\aaa \TrP \aab)}, &   & \aab = \aac,\\ [2mm]
      \aaa \TrP \aac,  &  & \aac \succ \aab,
    \end{array}\right.
\end{array}
$$
and compare the evaluations with respect to the different ordering
of the involved arguments. When elements of both $\Real$ and
$\tUnit$ are involved, use the above specification together with
Axiom \ref{ax:tropOperations}; for example,
$$
\begin{array}{llll}
 \uuu{\aaa} \TrP (\aab \TrS \aac)  =  & \uuu{(\aaa \TrP (\aab \TrS \aac))}
  & = &  \\[2mm]
   & \uuu{\OP(\aaa \TrP \aab) \TrS (\aaa \TrP \aac)\CP} & = &  \\[2mm]
   &  \uuu{\OP\aaa \TrP \aab\CP} \TrS \uuu{\OP\aaa \TrP \aac\CP} &=
   &
  (\aaa \TrP \uuu{ \aab} ) \TrS (\aaa \TrP \uuu{\aac}) \ . \\
\end{array}
$$
\parSpc
\noindent \bfem{Zero:} By definition $\zero_{\Trop} := \tUniS$ is
the additive identity of $\Trop$ (cf.~Axiom
\ref{ax:tropOperations} (1)), and it annihilates $\Trop$
(cf.~Axiom \ref{ax:tropOperations} (4)).

\parSpc
\noindent \bfem{One:} One can easily check that $\one_{\Trop} :=0$
is the multiplicative identity of $\Trop$.

\begin{theorem}\label{thm:propOfOperation} The set $\Trop$ equipped
with the addition $\TrS$ and the multiplication $\TrP$ is a
(non-idempotent) commutative semiring,  $(\tUnit, \TrS)$ is an
additive semigroup, and $(\Real, \TrP)$ and $(\tUnit, \TrP)$ are
multiplicative semigroups.
\end{theorem}

\begin{remark}\label{rmk:nu} In the view of Axiom \ref{ax:tropOperations}, $\nu$
is realized as the onto order preserving projection
\begin{equation}\label{eq:nu}
    \nu : (\Trop,\TrS,\TrP) \ \To \  (\etUnit,\TrS,\TrP) \ ,
\end{equation}
where $\nu : a \mapsto \uuu{a}$, $\homToU: \uuu{\aaa} \mapsto
\uuu{\aaa}$, and $\homToU:\tUniS \mapsto \tUniS$. Then, $\nu$ is a
semiring homomorphism and we write $\uuu{x}$ for the image of $x
\in \Trop$ in $\etUnit$, where $\homToU$ is  is the identity for
each $x \in \etUnit$. Accordingly, call $\uuu{a}$ the
\textbf{$\nu$-value} of $a$. Given $x,y \in \Trop$, we say that
$x$ is greater than $y$, or maximal, \textbf{up to} $\nu$ if
$\uuu{x} \succ \uuu{y}$, similarly, when $\uuu{x} = \uuu{y}$ we
say that $x$ and $y$ are equal up to $\nu$ .
\end{remark}
Writing  $x^n$ for the tropical product $x \TrP x \TrP \cdots \TrP
x$
 of $n$ factors we have:

\begin{lemma}
$
    (x \TrS y)^n = x^n \TrS y^n, \  n \in \Net,
$
for any $x,y \in \Trop$.
\end{lemma}
\begin{proof}
Assume $n > 1$, by  induction:
$$(x \TrS y)^n = (x \TrS y)(x \TrS y)^{n-1}=(x \TrS y)(x^{n-1} \TrS y^{n-1}) =
x^n \TrS x^{n-1}y \TrS xy^{n-1} \TrS y^{n}.$$
Suppose $x \succ y$, then
$$x^n \ \succ   \ x^{n-1}y  \ \TrS  \ xy^{n-1} \TrS
y^{n} \ \succ  \ y^{n}$$
and $(x \TrS y)^n  = x^n$. Similarly, if $y \succ x$, then $(x
\TrS y)^n  = y^n$. In the case of $x = y$, we have $x \TrS y \in
\etUnit$, $x^n \TrS y^n \in \etUnit$, and $x^n \TrS y^n = x^n = (x
\TrS y)^n$.
\end{proof}

\begin{corollary}\label{thm:powerOfPoly}
$\OP \bigoplus_{i=1}^s x_i \CP ^n = \bigoplus_{i=1}^s x_i^n, \
 n \in \Net,$ for any $x_1, \dots, x_s \in \Trop$.
\end{corollary}

\begin{corollary}\label{thm:CauchyInequality}
The ``Cauchy'' inequality
$$x_1 \TrP x_2 \TrP \cdots \TrP x_n \ \preceq \  x_1^n \TrS x_2^n \TrS \cdots \TrS x_n^n$$
holds for any $x_1, \dots, x_n \in \Trop$; equality occurs only if
$\homToU(x_1) = \homToU(x_2) = \cdots = \homToU(x_n)$ and at least
one $x_i$ is in $\etUnit$.
\end{corollary}

\subsection{Tropical arithmetics and
tropicalization}\label{sec:Tropicalization} The informal term
tropicalization is used to describe a map, based on a real
valuation, of objects defined over a non-Archimedean field $\Fld$
with real valuation to objects defined over $(\eReal, \max, + \
)$; objects are either varieties or polynomials. The
tropicalization of a variety $W \subset \Fld^{(n)}$ is a
polyhedral complex in $\Real^{(n)}$, while a polynomial in $n$
variables in $\Fld[\lm_1, \dots, \lm_n ]$ is mapped to a tropical
polynomial in $n$ variables in $\eReal[\lm_1, \dots, \lm_n ]$,
which we recall determines an affine piecewise linear function.

Let $\Fld$ be an algebraically closed field with a real
non-Archimedean valuation
\begin{equation}\label{eq:Valuation}
    \nVal : \; (\Fld,+,\cdot \ ) \  \To \  (\eReal,\max,+) \ ;
\end{equation}
for example,  assume $\Fld$ is the field of locally convergent
complex \bfem{Puiseux series}, of the form
$$
  f(t) = \sum_{\aaa \in R} \aac_{\aaa} t ^{\aaa},  \qquad \aac_{\aaa} \in
  \Comp \ ,
$$
where $R \subset \Rati$ is bounded from below and the elements of
$R$ have a bounded denominator. Then,
\begin{equation}\label{eq:valPowerSeries}
  \nVal(f) = \left\{%
\begin{array}{ll}
     - min \{\aaa \in R \ : \; \aac_{\aaa} \neq 0 \}, &  f \in  \Fld[\lm_1, \dots, \lm_n ] \setminus 0 \ ;
     \\ [2mm]
     \tUniS, &  f = 0 \ , \\
\end{array}%
\right.
\end{equation}
 is a real valuation satisfying the rules of being non-Archimedean,
\begin{equation}\label{eq:valRules}
\begin{array}{ll}
  (i)   &   \nVal(f \cdot g)  =  \nVal(f) + \nVal(g) \ , \\  [2mm]
  (ii) & \nVal(f + g)  \leq  \max \{ \nVal(f), \nVal(g)\} \ . \\
\end{array}
\end{equation}
(Note that $\nVal$ is not a homomorphism, since it does not
 preserve associativity.)
Thus, in the sense of tropicalization, the arithmetic operations
of $\Fld$ are replaced with the correspondence: $ \cdot \  \mapsto
\ + \ $ and $\ + \ \mapsto \ \max$ .

\begin{remark}
Taking $f,g \in \Fld$ with $\nVal(f) = \nVal(g) = \aaa$, then
$\nVal(f + g)$ can be any point of the ray $[\tUniS,\aaa]$.  These
cases provide the motivation for the use of $\tSR$ as the target
of $\nVal$ that allows to distinguish between the cases in which
Formula \Ref{eq:valRules}(ii) is interpreted as equality and the
cases it is inequality.
\end{remark}

In order to realize $\tSR$ as the target of $\Val$,  to each point
$\uuu{\aaa} \in \tUnit$ we assign the ray $P_{\uuu{\aaa}} : =
[\tUniS,\aaa]$ and to each $x \in \eReal$ we assign the singleton
$P_{{a}} := \{ a \}$, in particular $P_{{\tUniS}} := \{ \tUniS
\}$; therefore $x \in P_x$ for each  $x \in \Trop$. With this
construction we obtain the inclusions:
\begin{equation}\label{eq:pContainements}
  P_{\tUniS}, P_a \subset P_{\uuu{\aaa}}, \qquad
  P_{\uuu{a}}  \subset P_{\uuu{b}}  \Longleftrightarrow a \ \prec \ b,
  \qquad \forall a,b \in \Real \ .
\end{equation}
(Recall that two series in $\Fld$ that are vanished in order 1
must vanished on order at least 1; the inclusions
\Ref{eq:pContainements} address this property.)

Let $G(\Trop) := \{P_x \; : \; x \in \Trop \}$, then $\nVal(f) \in
P_x$ for some $P_x \in G(\Trop)$ which clearly needs not be
unique. Accordingly, for each pair $f \in \Fld$ and $x \in \Trop$
we define the relation
\begin{equation}\label{eq:homomorphicRelation}
  \nVal(f) \in P_x \qquad or \qquad \nVal(f) \notin P_x
\end{equation}
determined by the inclusion of $\nVal(f)$ in $P_x$.

\begin{theorem}\label{thm:homomorphicRelation}
Formula \Ref{eq:homomorphicRelation} yields a homomorphism; that
is, for any $f,g \in \Fld$ with $\nVal(f) \in P_x$ and $\nVal(g)
\in P_y$ we have $\nVal(f \cdot g) \in P_{x \TrP y}$ and $\nVal(f
+ g) \in P_{x \TrS y}$.
\end{theorem}
\begin{proof}
Suppose $ \nVal(f) = \aaa$, $\nVal(g) = \aab$. Then, since $x \in
P_x$ for each $x \in \Trop$,
$$
  \nVal(f \cdot g)  =  \nVal(f) + \nVal(g) =
   \aaa \TrP \aab  \in  P_{\aaa \TrP \aab} \ ,$$ and
   $
  \nVal(f \cdot 0)  =   \nVal(f) + \nVal(0) = \aaa \TrP (\tUniS)  \in
  P_{\tUniS}$ .
For the additive relation, write
$$
\begin{array}{lrrl}
  \nVal(f + g) \leq  & \max \{ \nVal(f), \nVal(g)\} & = &  \\
  & \max \{ \aaa, \aab \} & =  &
 \left\{
\begin{array}{lllllll}
  \aaa &\in& P_{\aaa} &=& P_{\aaa \TrS \aab} & & \aaa > \aab,  \\
  \aaa &\in& P_{\aaa} \subset P_{\uuu{\aaa}} & = & P_{\aaa \TrS \aaa}& & \aaa  =  \aab,\\
  \aab &\in& P_{\aab} &=& P_{\aaa \TrS \aab}  & & \aab > \aaa,
\end{array}
\right.
   \\
\end{array}
$$
 and use the inclusion  $P_{\aaa} \subset  P_{\uuu{\aaa}}$, cf. \Ref{eq:pContainements}. The case of
$\nVal(f + 0)$ is trivial.
\end{proof}
\subsection{The relation to the max-plus arithmetic}
\label{sec:ExtendedArithmetic} The structure of $\tSR$ provides a
much richer structure,  generalizing both the max-plus semiring
and the one  suggested in Remark \ref{rem:nonAssociative}, and
achieves the best of both worlds.

\begin{lemma}\label{lem:homomorphismsMax2Extended}
The map
\begin{equation}\label{eq:homomorphismsMax2Extended}
    \epiToMaxPlus: (\Trop,\TrS,\TrP) \  \To \ (\eReal, \max,+ \ ) \ , \\
\end{equation}
%
$\epiToMaxPlus:\uuu{\aaa} \mapsto \aaa$, $\epiToMaxPlus:\aaa
\mapsto \aaa$, and $\epiToMaxPlus:\tUniS \mapsto \tUniS$, is a
semiring epimorphism.
\end{lemma}
\begin{proof} Clearly, $\epiToMaxPlus$ is onto. Assume $\epiToMaxPlus(x) =
a$ and $\epiToMaxPlus(y) = b$, where $x, y \in \Trop$, then
$\epiToMaxPlus(x  \TrS y) = \max\{a,b\} = \max\{\epiToMaxPlus(x),
\epiToMaxPlus(y) \}$ and $\epiToMaxPlus(x \TrP y) = a +b =
\epiToMaxPlus(x) + \epiToMaxPlus(y) $.
\end{proof}

On the other hand one can also define:
\begin{lemma}\label{thm:isoRelaToU}
The  map
\begin{equation}\label{eq:maxPlusToU}
    \maxPlusToU: (\eReal, \max,+ \ ) \  \To  \ (\etUnit, \TrS,\TrP), \\
\end{equation}
%
 $\maxPlusToU:\aaa \mapsto
\uuu{\aaa}$ and $\maxPlusToU:\tUniS \mapsto \tUniS$, is a semiring
isomorphism that embeds $(\eReal, \max,+ \ )$ in $(\Trop,
\TrS,\TrP)$.
\end{lemma}
\begin{proof}
Take  $\aaa, \aab \in \eReal$, then $\maxPlusToU( max\{ \aaa,
\aab\})  = \uuu{(\max\{ \aaa, \aab\})} = \uuu{\aaa} \TrS
\uuu{\aab} =  \maxPlusToU(\aaa) \TrS \maxPlusToU(\aab),$ and $
\maxPlusToU( \aaa + \aab)  = \uuu{(\aaa +\aab)} = \uuu{\aaa} \TrP
\uuu{\aab} =  \maxPlusToU(\aaa) \TrP \maxPlusToU(\aab).$ $\etUnit
\subset \Trop$, so $\maxPlusToU$ embeds  $(\eReal, \max,+ \ )$ in
$(\Trop, \TrS,\TrP)$.
\end{proof}

\begin{corollary}\label{cor:gen}
Categorically, by Remark \ref{rmk:nu}, Lemma
\ref{lem:homomorphismsMax2Extended},  and Lemma
\ref{thm:isoRelaToU}, the diagram
$${\setlength{\unitlength}{0.7mm}
\begin{picture}(60, 27)
  \put(10,22){$(\Trop,\TrS, \TrP)$}
  \put(55,22){$(\eReal,\max, + \ )$}
  \put(57,0){$(\etUnit,\TrS, \TrP  )$}

  \put(40,24){$\epiToMaxPlus$}
  \put(40,12){$\nu$}
  \put(69,12){$\maxPlusToU$}

  \put(32, 23){\vector(1, 0){20}}
  \put(20, 20){\vector(2, -1){35}}
  \put(67, 20){\vector(0, -1){15}}
\end{picture}}$$
commutes.
\end{corollary}

Corollary \ref{cor:gen} displays  $\tSR$ as a generalization of
$(\eReal,\max, + \ )$ which is endowed with a richer structure in
the sense that it encodes an indication about the additive
multiplicity of elements in $\Real$. Namely, since $a \TrS a =
\uuu{a}$ and $a \TrS \uuu{a} = \uuu{a}$, $\uuu{a}$ can be realized
as a point with additive multiplicity $>1$. (Clearly, computations
for $(\eReal, \max,+ \ )$ can be performed on $\tSR$ and then to
be sent back to $(\eReal, \max,+ \ )$.)

As for the  arithmetic suggested in Remark
\ref{rem:nonAssociative} (i.e. defined with $``\aaa + \aaa"  =
\tUniS$), one may suggest the map
\begin{equation}\label{eq:homomorphismsTrop2Extended}
    \phi: (\Trop,\TrS,\TrP) \ \To \ (\eReal, ``\max",+  \ ), \\
\end{equation}
$\phi:\aaa \mapsto \aaa$, $\phi:\uuu{\aaa} \mapsto \tUniS$, and
$\phi:\tUniS \mapsto \tUniS$; but, since $(\eReal, ``\max",+ \ )$
is not associative, $\phi$ is not a homomorphism.


\subsection{Geometric view} \label{sec:GeometricView}

Let us remind  that one of our goals was to obtain a semiring
structure that enables us to treat algebraically the points of a
corner locus of tropical functions, namely,  to define tropical
algebraic set. To present only the frame of this idea, given a
tropical polynomial $f\in \Trop[\lm_1,\dots,\lm_n]$ we define the
tropical algebraic set  of the corresponding function $\tilde{f}$
to be
$$Z(\tilde f) = \{ (x_1,\dots,x_n) \in \Trop^{(n)} \ : \ \tilde f(x_1,\dots,x_n) \in \etUnit \}.$$
Therefore, the corner locus of $f \in \eReal[\lm_1,\dots,\lm_n]$
over $(\eReal, \max, + \ )$ is just the restriction of $Z(\tilde
f)$,  considered as a polynomial over $(\Trop, \TrS, \TrP )$, to
the  real points, i.e. $Z(\tilde f) \cap \eReal^{(n)}$.
\begin{example}

\begin{figure}[!h]
\setlength{\unitlength}{0.7cm}
\begin{picture}(18,8)(0,0)

\thinlines

\drawline[-20](0,4)(8,4) \drawline[-20](11,4)(17,4)
\drawline[-20](4,1)(4,7) \drawline[-20](14,1)(14,7)

\thicklines \drawline(11.5,5)(15,5)(17,7)
\dottedline{0.15}(15,5)(15,4)
\drawline(5.6,5.1)(7.5,7)\drawline(5.4,5.05)(4.8,5.05)\drawline(2.7,5.05)(3.3,5.05)
\drawline(2.6,2.9)(1,1.3)\dottedline{0.15}(4.8,5.05)(3.3,5.05)

\put(17.2,3.9){$x \in \Real$} \put(8.2,3.9){$x \in \Trop$}
\put(13.8,7.3){$\tilde f(x)$}\put(3.8,7.3){$\tilde
f(x)$}\put(14.9,3.6){$a$}
\put(5.4,4.9){$\circ$}\put(2.5,4.9){$\circ$}
\put(2.5,2.8){$\bullet$}\put(5.4,2.8){$\bullet$}
\put(5.3,3.6){$a$}\put(2.4,3.6){$a^\nu$}\put(3.3,3.6){$-\infty$}

\put(3.3,0.2){$(\Trop, \TrS, \TrP)$} \put(13.3,0.2){$(\eReal,
\max, + \ )$}
\end{picture}
\caption{\label{fig:basicGraph} The graph of the linear function
$\tilde f(x) = x + a $ over $(\Trop, \TrS,\TrP)$, on left hand
side, and the corresponding function $\tilde f(x) = \max\{x,a \}$
over $(\eReal, \max,+ \ )$, on right hand side.}
\end{figure}
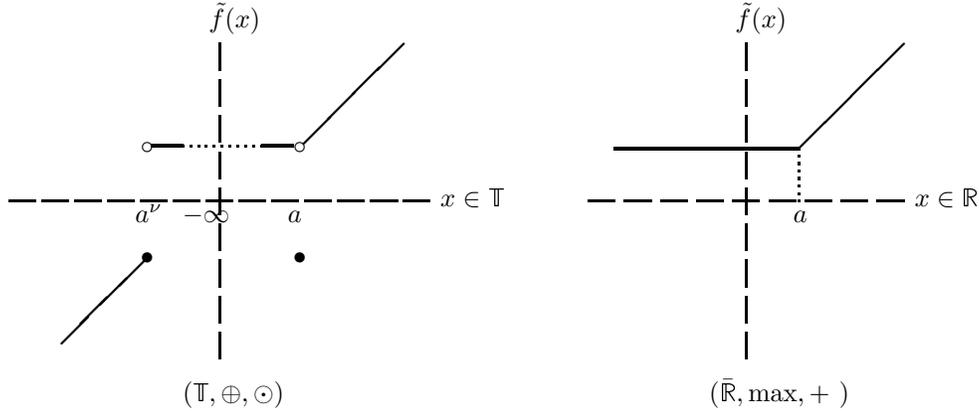
Consider the similar linear functions $f(x) = x \TrS \aaa$ over
$(\Trop,\TrS,\TrP)$ and $f(x) = \max\{ x, \aaa\}$ over $(\eReal,
\max,+)$, see Figure \ref{fig:basicGraph}. Restricting the domain
to $\Real$ only, over  $(\Trop, \TrS,\TrP)$, the image of the
corner locus, which contains the single point $\aaa$, is
distinguished and is now mapped to $\tUnit$.
\end{example}
The study of polynomial algebras and tropical algebraic sets over
$(\Trop,\TrS,\TrP)$ will be treated in a forthcoming paper.

\section{Matrix Algebra}\label{chap:TropicalAlgebra}

Our forthcoming study is dedicated to introducing the fundamentals
of the matrix algebra over $\tSR$ whose operations of  are
typically combinatorial.  Yet, developing an
 algebraic theory, analogous to classical theory of matrix algebra over fields,
 with a view to combinatorics, is our main
 goal. This goal is supported by the connections to graph theory \cite{Lawler76},
 the theory of automata
 \cite{margolisEdense}, and semiring theory \cite{golan92}.
\pSkip

\bfem{Notations:} For the rest of the paper, assuming the nuances
of the different arithmetics are already familiar, we write $xy$
for the product $x \TrP y$,   $\frac{\ x \ }{ \ y \ }$ for the
division $x \iTrP y$, and $x^n$ for $x \TrP \cdots \TrP x$
repeated $n$ times.


\subsection{Tropical matrices}\label{sec:FundamentalofTropicalMatrices}

It is standard that if $R$ is a semiring then we have the semiring
$M_n(R)$ of $n \times n$ matrices with entries in $R$, where
addition and multiplication are induced from $R$ as in the
familiar matrix construction. Accordingly, we define the semiring
 of tropical matrices $M_n(\Trop)$ over $(\Trop, \TrS,
\TrP)$, whose unit is the matrix
\begin{equation}\label{def:absUnitMat}  
\tI = \OP  \begin{array}{ccc}
  0 & \ldots & \tUniS \\
  \vdots & \ddots & \vdots \\
  \tUniS & \ldots & 0
\end{array} \CP \\
\end{equation}
and whose zero matrix is  $\tZ = (\tUniS) \tI$; therefore,
$M_n(\Trop)$ is also a multiplicative monied. We write  $\tA
=(a_{ij})$ for a tropical matrices $\tA \in M_n(\Trop)$  and
denote the entries of $\tA$ as $a_{ij}$. Since $\Trop$ is a
commutative semiring, $ x \tA = \tA x$ for any $x \in \Trop$ and
$\tA \in M_n(\Trop)$.

As in the familiar  way, we define the \bfem{transpose} of $\tA =
(a_{ij})$ to be $\tA^\trn = (a_{ji})$, and have the relation
\begin{proposition}\label{thm:commutativeOfSimmMat}
$\OP \tA \tB \CP  {^\trn} =   \tB^\trn   \tA^\trn. $
\end{proposition}
\noindent  (The proof is standard by the commutativity and the
associativity of $\TrS$ and $\TrP$ over $\Trop$.)

The \bfem{minor} $\tA_{ij}$ is obtained by deleting the $i$ row
and $j$ column of $\tA$. We define the \bfem{tropical determinant}
to be
  \begin{equation}\label{def:tropicalDet}
 |\tA| = \bigoplus_{\sig \in S_n}  \OP  \aaa_{1\sig(1)}
 \cdots \aaa_{n\sig(n)} \CP,
\end{equation}
where $S_n$ is the set of all the permutations on $\{1,\dots,
  n\}$.   Equivalently,  $|\tA|$ can be written in terms of minors as
\begin{equation}\label{def:tropicalDetByMinors}
|\tA| = \bigoplus_{j}   \aaa_{i_o j} |\tA_{i_o j}|,
\end{equation}
for some fixed index $i_o$. Indeed, in the classical sense, since
parity of indices' sums are not involved in Formula
\Ref{def:tropicalDet}, the tropical determinant is a permanent,
which makes the tropical determinant a pure combinatorial
function. The \bfem{adjoint} matrix $\Adj{\tA}$ of $\tA =(a_{ij})$
is defined as the matrix $(a'_{ij})^{\operatorname{t}}$ where
$a'_{ij}= |\tA_{ij}|$.

\begin{observation}\label{obs:detProp}
The tropical determinant has the following properties:
\begin{enumerate}
   \item Transposition and reordering of rows or columns leave the determinant
   unchanged; \pSkip
  \item The determinant is linear with respect to scalar multiplication of any given row or column.
\end{enumerate}
\end{observation}

\subsection{Regularity of matrices } Using the special structure
of $(\Trop, \TrS, \TrP)$, the algebraic formulation of
combinatorial properties becomes possible.

  \begin{definition}\label{def:singularMat} A matrix $\tA \in M_n(\Trop)$
  is said to be \textbf{tropically singular}, or singular, for short,
  whenever
  $|\tA| \in \etUnit$, otherwise $\tA$ is called \textbf{tropically regular},
  or regular, for short.
  \end{definition}
  %
 In particular, when two or more different
permutations, $\sig \in S_n$, achieve the $\nu$-value of $|\tA|$
simultaneously, or the permutation that reaches the $\nu$-value of
$|\tA|$ involves an entry in $\etUnit$, then $\tA$ is singular.

\begin{remark}\label{obs:regProp}
Despite some classical properties hold for the tropical
determinant, cf. Observation~\ref{obs:detProp},  the familiar
relation $|\tA \tB| = |\tA| |\tB|$ does not hold true on our
setting; for example, take the matrix
  \begin{equation}\label{exmp:DetMultiplication}
   \tA = \vMat{1}{1}{2}{3}  \quad \text{ with } \quad
  \tA^2 = \tA  \tA =
  \vMat{3}{4}{5}{6},
\end{equation}
   then,  $|\tA| = 4$ and $|\tA| |\tA| =
   8$, while $|\tA^2| = \uuu{9}$. In the view of tropicalization, which ignores signs,
   the determinant of a matrix over $\Fld$ is assigned to the permanent of a matrix
   in $M_n(\Trop)$;
   this explains the tropical situation in which the product of two regular matrices might be singular.
 \end{remark}

\begin{theorem}\label{thm:detof2rows}
A matrix with two identical rows or columns is singular.
\end{theorem}
\begin{proof} Proof by induction on $n \geq 2$. The case of $n=2$ is clear
by direct computation. Assume the two first columns of $\tA$ are
identical, and expand $|\tA|$ in terms of minors along the first
row, that is $|\tA| = \bigoplus_i a_{1i}|\tA_{1i}|$. Since $a_{11}
= a_{12}$ and $\tA_{11} = \tA_{12}$, then $a_{11}|\tA_{11}| =
a_{12}|\tA_{12}|$, and so  $a_{11}|\tA_{11}| \TrS a_{12}|\tA_{12}|
\in \etUnit$. By the induction hypothesis, for any $i >2$,
$\tA_{1i}$ is a matrix with identical columns, and is singular,
that is  $|\tA_{1i}| \in \etUnit$ for all $i > 2$. When adding all
together, $a_{1i}|\tA_{1i}| \in \etUnit$ for all $i = 1,\dots,n$,
and thus $|\tA| \in \tUnit$.
\end{proof}

\begin{theorem}\label{thm:detofMult}
If $\tA$ and $\tB$ are regular matrices and their product $\tA\tB$
is also regular, then $ |\tA\tB| = |\tA| |\tB|$. When either $\tA$
or $\tB$ is singular, then $\tA\tB$ is also singular.

\end{theorem}

\begin{proof}
Let $S_n$ be the set of all the permutations on $N =
\{1,\dots,n\}$, and let $\maps_n = \{N \to N\}$ be the set of all
maps from $N$ to itself, in particular $S_n \subset \maps_n$.
Denoting the entries of $\tA \tB$ by $(ab)_{ij}$, we write the
determinant $|\tA \tB|$ in the form of Formula
\Ref{def:tropicalDet} as:
$$|\tA\tB| = \bigoplus_{\sig \in S_n} \bigodot_{i} (\aaa \aab)_{i\sig(i)} =
 \bigoplus_{\sig \in S_n} \bigodot_{i}\OP \bigoplus_k (\aaa_{ik} \aab_{k\sig(i)})\CP = $$
$$
\bigoplus_{\sig \in S_n} \OP
   (\aaa_{11}  \aab_{1\sig(1)} \TrS  \cdots \TrS \aaa_{1n}
   \aab_{n\sig(1)}) \
   \cdots \
   (\aaa_{n1}  \aab_{1\sig(n)} \TrS  \cdots \TrS \aaa_{nn}  \aab_{n\sig(n)})
 \CP =
$$
 \begin{equation}\renewcommand{\theequation}{*}\addtocounter{equation}{-1}\label{eq:str.1}
\bigoplus_{\sig \in S_n}   \bigoplus_{\mu \in \maps_n} \OP
  \bigodot_i \OP \aaa_{i\mu(i)}  \aab_{\mu(i)\sig(i)} \CP \CP  =
 \bigoplus_{\sig \in S_n}   \bigoplus_{\mu \in \maps_n}
 \OP \bigodot_i \aaa_{i\mu(i)} \bigodot_i \aab_{\mu(i)\sig(i) }
 \CP.
\end{equation}
By the structure of the left hand side of $(*)$, we can see that
the value of $|\tA \tB|$ is obtained when both $\bigodot_i
\aaa_{i\mu(i)}$ and $\bigodot_i \aab_{\mu(i)\sig(i)}$ attain their
maximal evaluation  at the same time. We show that this is
possible. Namely both reach their maximal evaluation on the same
$\mu$, which we denote by $\mu_o$; the corresponding $\sig \in
S_n$ is then denoted by $\sig_o$.  Note that when $|\tA\tB| \in
\Real$ there must be exactly one pair, $\mu_o$ and $\sig_o$;
otherwise, by definition, $\tA\tB$ would not be regular.

\parSpc
 \noindent \bfem{Case I:}  Suppose  $\mu_o \in S_n$ is a
 permutation which
  maximizes $\bigodot_i \aaa_{i\mu(i)}$. We
 show that there is also a permutation $\sig_o \in S_n$ that maximizes
$\bigodot_i \aab_{\mu(i)\sig(i)}$ for $\mu_o$. Assume $|\tA\tB|
\in \Real$, with $\sig_t \in S_n$ maximizes $\bigodot_j
\aab_{j\sig_t(j)}$. Generally speaking, for any given $\mu \in
S_n$ and $\sig_t \in S_n$, there exits $\sig \in S_n$ which makes
the diagram
\begin{equation*}
  \begin{CD}
    \qtriangle[N_{[j]}`N_{[i]}`N;\mu`\sig_t`\sig]
  \end{CD}
\end{equation*}
commutative, where we use the notation ${[ \; \cdot \; ]}$ to
indicate the appropriate  indices. Accordingly, choosing $\sig_o
\in S_n$ for which  $\sig_o \circ \mu_o = \sig_t $, we obtain
$\bigodot_i \aab_{\mu_o(i)\sig_o(i)} = \bigodot_j
\aab_{j\sig_t(j)}$; in this case the two components of $(*)$ reach
their maximum simultaneously and we can write:
$$ 
  (*) =
   \OP \bigodot_i \aaa_{i\mu_o(i)} \CP
   \OP \bigodot_j \aab_{j\nu_o(j)} \CP =
  \OP \bigoplus_{\mu \in S_n} \bigodot_i \aaa_{i\mu(i)} \CP
  \OP \bigoplus_{\nu \in S_n} \bigodot_j \aab_{j\nu(j)} \CP =
  |\tA| |\tB| \ .
$$

When $\tA$ is singular, there are at least
 two different $\mu_1,\mu_2 \in S_n$ that attain the $\nu$-value of $|\tA|$
 or a
 single  $\mu \in S_n$ that involves a non-real entry.
The latter case is obvious, since $(*)$ has a non-real multiplier
and thus $|\tA \tB| \in \etUnit$.
 Suppose  $\sig_1,\sig_2 \in S_n$ are two
permutations satisfying $  \nu_o  = \sig_l\circ \mu_l$, $l=1,2$,
then
$$ (*) =
  \bigodot_i \aaa_{i\mu_1(i)} \bigodot_i \aab_{\mu_1(i)\sig_1(i) } =
  \bigodot_i \aaa_{i\mu_2(i)} \bigodot_i \aab_{\mu_2(i)\sig_2(i) } \ ,
$$
 and hence
$|\tA\tB| \in \etUnit$.

\noindent
  \bfem{Case II:}
 Suppose $\mu_o \in \maps_n \setminus S_n$, $|\tA\tB| \in \Real$, and let $\sig_o \in S_n$ be the
  corresponding permutation which maximizes the product
  \begin{equation}\renewcommand{\theequation}{**}\addtocounter{equation}{-1}\label{eq:str.2}
\bigodot_i \OP \aaa_{i\mu(i)} \aab_{\mu(i)\sig(i) } \CP =
     \bigodot_i \aaa_{i\mu(i)} \bigodot_i \aab_{\mu(i)\sig(i) } \ .
\end{equation}
in $(*)$.   In particular, there is only one such pair, $\mu_o$
and $\sig_o$, for
  otherwise $\tA\tB$ would not be regular.
  Since $\mu_o \notin S_n$, there are at least two indices $i_1
  \neq i_2$ with
$ \mu_o(i_1) = \mu_o(i_2) = k_o$. Let $h_1 := \sig_o(i_1)$ and
$h_2 := \sig_o(i_2)$; then $h_1 \neq h_2$, since $\sig_o \in S_n$.
Subject to $\mu_o$ and $\sig_o$, $(**)$ can be rewritten as,
   $$ \OP \bigodot_i \aaa_{i\mu_o(i)}       \CP
      \OP \bigodot_i \aab_{\mu_o(i)\sig_o(i)}\CP =
     \OP \bigodot_{i}  \aaa_{i\mu_o(i)} \CP
     \OP (\aab_{\mu_o(i_1),\sig_o(i_1) } \aab_{\mu_o(i_2)\sig_o(i_2) })
     \bigodot_{i \neq i_1 i_2} \aab_{\mu_o(i)\sig_o(i) }
     \CP =
  $$

 \begin{equation}\renewcommand{\theequation}{***}\addtocounter{equation}{-1}\label{eq:str.1}
    \OP \bigodot_{i}  \aaa_{i,\mu_o(i)} \CP
     \OP (\aab_{k_o h_1} \aab_{k_o h_2}) \bigodot_{i \neq i_1 i_2} \aab_{\mu_o(i) \sig_o(i) }
     \CP.
\end{equation}
 Denote by $\tilde{\sig}_o \in S_n$ the permutation obtained
  by switching between the images of $i_1$
and $i_2$ in $\sig_o$, while all other correspondences remain as
they are;  explicitly, $\tilde{\sig}_o(i_1) = h_2$,
$\tilde{\sig}_o(i_2) = h_1$, and $\tilde{\sig}_o(i) = \sig_o(i)$,
for all $i \neq i_1, i_2$. With respect to $\sig_o$, we have the
following combinatorial situation:
$$
\begin{array}{|cc|c|ccc|c|c|}
\hline
     &  & \overset{\sig_o(i_1)}{=h_1} &  &  &  & \overset{\sig_o(i_2)}{=h_2}  & \qquad \\
     &  &  &  &  &  &  & \\
     &  & \vdots &  &  &  & \vdots & \\\hline
    \begin{array}{c}
       \overset{\mu_o(i_1)=}{} \\
       \overset{\mu_o(i_2)= }{} \
     \end{array}  k_o
     & \cdots & * &  &  \cdots&  &  *  & \cdots \\\hline
     &  & \vdots &  &  &  & \vdots & \\
     &  &  &  &  &  &  & \\
     &  & \overset{=\tilde{\sig}_o(i_2)}{} &  &  &  & \overset{= \tilde{\sig}_o(i_1)}{} & \\ \hline
\end{array}$$
(The diagram helps us to understand the modification of $\mu_o$.)
 Using $\tilde{\sig}$ we  expand  $(***)$ further,
   $$(***) =
     \OP \bigodot_{i}  \aaa_{i\mu_o(i)} \CP
     \OP (\aab_{\mu_o(i_2)\tilde{\sig}_o(i_2) }  \aab_{\mu_o(i_1)\tilde{\sig}_o(i_1)
     })\bigodot_{i \neq i_1 i_2} \aab_{\mu_o(i)\tilde{\sig}_o(i) }
     \CP =
  $$
  $$
     \OP \bigodot_{i}  \aaa_{i\mu_o(i)} \CP
     \OP \bigodot_{i } \aab_{\mu_o(i)\tilde{\sig}_o(i) } \CP,
  $$
which means that $\tilde{\sig}_o \neq \sig_o$ also attain the
$\nu$-value of $|\tA \tB|$ and thus,  $|\tA\tB| \in \etUnit$. This
contradicts the assumption that $\mu_o \notin S_n$, so $\mu_o \in
S_n$, and this case has already been discussed before.
\end{proof}

\begin{example} Take the matrices
$$
    \tA = \vMat{1}{2}{2}{3} \quad \text{ and }
  \quad \tB = \vMat{3}{1}{0}{2},  \quad \text{then} \quad \tA \tB =
  \vMat{4}{4}{5}{5} \ .
$$
 $\tA$ is singular
 with $|\tA| = \uuu{4}$, $\tB$ is regular with $|\tB| = 5$, and
  $\tA \tB$ is singular with $ |\tA\tB| = \uuu{9}$. On the
 other hand, $\tB^2 =  \vMat{6}{4}{3}{4} $ is regular with $|\tB^2| =
 10$, so $|\tB| |\tB| = |\tB^2|$.

\end{example}
\section{Invertibility of Matrices}\label{sec:TheTropicalInverseMatrix}
We introduce a new notion of semigroup invertibility, and present
it for the matrix monoid $M_n(\Trop)$; this type of invertibility
can be adopted to any abstract semigroup having a distinguished
subset. Although our framework is typically combinatorial, we show
how classical results are carried naturally on our setting.

\subsection{Basic definitions}\label{sec:DefinitionandConstruction}
We open with an abstract definition.
\begin{definition}\label{def:semigInve}
Let $\tS$ be semigroup, and let $U \subset S$ be a proper subset
with the property that for any $u \in U$ there exists some $v\in
U$ for which $vu \in U$ and $uv \in U$. We call $U$ a
\textbf{distinguished subset} of $S$.

An element $x \in \tS$ is said to be \textbf{pseudo invariable} if
there is $y \in \tS$ for which $xy \in U$ and $yx \in U$, in
particular all the members of $U$ are pseudo invertible. When $U$
consists of all idempotents elements of $\tS$, the pseudo
invertibility is then called \textbf{E-denseness}
\cite{margolisEdense}. A monoid is called \textbf{E-dense} if all
of its elements are E-denseness.
\end{definition}
\noindent To emphasize, for the purpose of pseudo invertibility,
$U$ needs not  be closed under the law of $\tS$. The notion of
E-denseness is already known in literature, while the weaker
version of pseudo invertibility is new.

To apply the notion of pseudo invariability to $M_n(\Trop)$,
viewed as monoid, we define a \bfem{pseudo unit matrix} to be a
\textbf{regular} matrix of the form
\begin{equation}\label{eq:parUnitMat}
\gtI =  \OP  \begin{array}{ccc}
  0 & \ldots & \uuu{\io}_{ij} \\
  \vdots & \ddots & \vdots \\
  \uuu{\io}_{ji} & \ldots & 0
\end{array} \CP.
\end{equation}
that is $\io_{ij} \in \etUnit$ for all $i \neq j$, and $\io_{ii} =
0$, for each  $i = 1,\dots,n$. Since $\gtI$ is regular we
necessarily have $|\gtI| = 0$, and in particular the unit matrix
$\tI$, cf. \Ref{def:absUnitMat}, is also a pseudo unit. We define
the distinguished subset $U_n(\Trop) \subset M_n(\Trop)$ to be
\begin{equation}\label{eq:setOfPartialUnits}
  U_n(\Trop) = \left\{\gtI \; : \; \gtI \; \;  \text{ is a  pseudo unit
  matrix}
  \right\} \ ,
\end{equation}
Therefore, $\tI \in  U_n(\Trop)$ and  hence $\tI \gtI = \gtI\tI =
\gtI$, for each $\gtI \in  U_n(\Trop)$, which makes $U_n(\Trop)$ a
distinguished subset satisfies the condition of Definition
\ref{def:semigInve}.

Correspondingly, we define the distinguished subset
$U_n^{idem}(\Trop) \subset M_n(\Trop)$ to be
\begin{equation}\label{eq:setOfIdemPartialUnits}
  U_n^{\idem}(\Trop) = \left\{\gtI \; : \; \gtI \; \text{ is an idempotent  pseudo unit
  matrix}
  \right\}\ .
\end{equation}
\begin{remark} It easy to show that any $\gtI \in U_2(\Trop) $ is
idempotent. For $n > 2$, not all of the pseudo units are
idempotents; for example, take the triangular matrix $$\gtI =
\vvMat{0}{\uuu{a}}{\uuu{b}}{\tUniS}{0}{\uuu{c}}{\tUniS}{\tUniS}{0}
\ ,$$ with $\uuu{a}\uuu{c} \succ \uuu{b}$.
\end{remark}

Using $U_n(\Trop)$ we explicitly define  pseudo invertibility on
$M_n(\Trop)$:
\begin{definition}\label{def:inverseMatrix}
A matrix $\tA \in M_n(\Trop)$ is said to be pseudo invertible  if
there exits a matrix $\tB \in M_n(\Trop)$ such that $\tA \tB \in
U_n(\Trop)$ and $\tB \tA \in U_n(\Trop)$. If $\tA$ is pseudo
invertible, then we call $\tB$  a \textbf{pseudo inverse matrix}
of $\tA$ and denote it as $\itA$.
\end{definition}
  We use the  notation of $\itA$ since the pseudo
matrix needs not be unique; moreover, in our setting $\tA \itA $
is not necessarily equal to $\itA \tA$, and thus might be
evaluated for different pseudo units.
\begin{example} Consider the following matrices:
$$\small{\tA = \OP  \begin{array}{ccc}
  0  & -2 & -1 \\
  -2 &  0 & \uuu{(-3)} \\
  -1 & \uuu{(-3)} & 0
\end{array} \CP  \qquad and  \qquad \tA' =  \OP  \begin{array}{ccc}
  0 & -2& -1 \\
  -2 & 0 & -3 \\
  -1 & -3 & 0
\end{array} \CP.}
$$
For these matrices we have, $\tA \tA' \in U_n(\Trop)$, $\tA' \tA
\in U_n(\Trop)$, and also $\tA
 \tA \in U_n(\Trop)$. Namely, $\tA$ has at least two pseudo inverses.
\end{example}

\begin{remark} For the case of $M_2(\Trop)$,
our notion of pseudo invertibility coincides with the notion of
general invertibility in a semigroup in the sense of Von-Neumann
regularity \cite{Lallement}, but not for $M_n(\Trop)$ with $n>2$.
\end{remark}

\begin{remark}
When one intends to use the other semirings structures, either
$(\eReal, \max,+ \ )$ or $(\eReal, ``\max",+ \ )$, in order to
define an inverse matrix, or pseudo inverse, it appears to be very
restricted  or even impossible. Over $(\eReal, \max,+ \ )$, unless
$\tUniS$ is involved, the zero element $\tUniS$ is unreachable by
tropical sums and products of entries. Thus, obtaining the unit
 matrix $\tI$ as products of matrices is very restricted. On the other
hand, when using $(\eReal, ``\max",+ \ )$, as suggested in Remark
\ref{rem:nonAssociative}, multiplication is not associative, which
makes implementation very difficult.
\end{remark}
\subsection{Theorem on tropical pseudo inverse matrix}\label{sec:TheoremofTropicalInverseMatrix}

\begin{theorem}\label{thm:inverseMatrix}
A matrix $\tA \in M_{n }(\Trop)$ is pseudo invertible if and only
if is tropically regular. In case $\tA$ is regular, $\itA$  can be
defined as
$$\itA = \frac{Adj(\tA)}{|\tA|} \ .$$
\end{theorem}
Before proving the theorem, we recall some definitions and present
new notation:
\renewcommand{\labelenumi}{(\alph{enumi})}

\begin{enumerate}
\item  Division in $\Trop$ is denoted by $\frac{\;
\cdot
  \;}{\cdot}$ and interpreted as the substraction $\aaa -\aab$ in the classical
  sense. We write $a^{-1}$ for  $\frac{0}{a}$ and $a^m$ for the
  tropical product of $a$ repeated $m$ times (which is just $m \cdot a$ in the usual
  sense).
\pSkip
      \item   We use the notation $\tA_{ih,jk}$  for $(\tA_{ij})_{hk}$, that is
  $(h,k)$-minor of the minor $\tA_{ij}$, where  $ h \neq i$ and  $ k \neq j$
  with respect to the initial indices of $\tA$. Accordingly,  $|\tA_{ij}|$
  is written in terms of minors as
$|\tA_{ij}| =  \bigoplus_{k\neq j} a_{hk}  |\tA_{ih,jk}|$, where
$h \neq i. $

  \end{enumerate}
  \renewcommand{\labelenumi}{(\arabic{enumi})}

\begin{proof}
We prove only  multiplication on right,
$\tA \itA \in U_n(\Trop)$;
the multiplication on left is proved in the same way.

$(\Leftarrow)$ Assume $|\tA| \in \etUnit$ and at the same time
there exists a pseudo inverse $\itA$. Then, by Theorem
\ref{thm:detofMult}, $|\tA \itA| \in \etUnit$ and their product is
singular.  Recalling that $|\gtI| = 0$ for all $\gtI \in
U_n(\Trop)$, we have  $\tA \itA \notin U_n(\Trop)$.
\parSpc

$(\Rightarrow)$ We write $\stA$ for the adjoint matrix $\Adj{\tA}
$, for short, and denote the product $\tA \stA$ by $\tB =
(b_{ij})$. Assuming $\tA$ is regular, we need to  prove that
$ \frac{\tA\stA}{|\tA|} \in U_n(\Trop)$, or equivalently, that
$\tA \stA = |\tA|\gtI$
for some $\gtI \in U_n(\Trop)$.  To prove this, we need to verify
 the following conditions:
\begin{enumerate}
      \item $ b_{ii} = |\tA|$ for each $i$; \pSkip
    \item $b_{ij} \in \etUnit$, for any $i\neq j$;\pSkip
    \item $|\frac{\tB}{|\tA|}|  = 0$. \pSkip
\end{enumerate}

\noindent \bfem{Diagonal entries:} When $i=j$,
\begin{equation}\label{proof:1}
b_{ii} = \bigoplus_{k} \pa_{ik} \sia_{ki} =
    \bigoplus_{k} \pa_{ik} |\ptA_{ik}| = |\ptA|,
\end{equation}
since this is just the expansion of $|\tA|$ along row $i$ (cf.
Equation~\Ref{def:tropicalDetByMinors}).

%
%
%
\noindent
  \bfem{Non-diagonal entries:} For $i\neq j$,
   \begin{equation}\label{proof:2}
   b_{ij} = \bigoplus_{k} \pa_{ik} \sia_{kj} =
   \bigoplus_{k} \pa_{ik} |\ptA_{jk}|  \in \etUnit,
   \end{equation}
since this is the expansion of the determinant of the matrix
obtained from $\tA$ by replacing row $j$ with a copy of row $i$,
and which therefore has two identical row and is singular (Theorem
\ref{thm:detof2rows}).

\parSpc \noindent
\bfem{{Regularity of product}:} To prove  $|\frac{\tB}{|\tA|}| =
0$, we show equivalently that $|\tB| = |\tA|^n$. Let $S_n$ be the
set of all permutations on $N = \{1,\dots,n\}$ and let $\maps_n =
\{N \to N\}$ be the set of all maps from $N$ to itself, i.e. $S_n
\subset \maps_n$, and write the expansion of $|\tB|$ explicitly,
$$|\tB| =  \bigoplus_{\sig \in S_n} \OP \bigodot_i b_{i\sig(i)} \CP =
\bigoplus_{\sig \in S_n} \OP \bigodot_i \OP \bigoplus_{k} \pa_{ik}
|\ptA_{\sig(i)k}| \CP \CP = $$
$$\bigoplus_{\sig \in S_n} \OP
   \OP \pa_{11}  |\ptA_{\sig(1)1}| \TrS \cdots  \TrS \pa_{1n}  |\ptA_{\sig(1)n}| \CP
\  \cdots \
   \OP \pa_{n1}  |\ptA_{\sig(n)1}| \TrS \cdots  \TrS \pa_{nn}  |\ptA_{\sig(n)n}| \CP
\CP =$$
\begin{equation}\label{proof:3}
 \bigoplus_{\sig \in S_n}\bigoplus_{\mu \in \maps_n}
 \bigodot_i \OP  \pa_{i\mu(i)} |\ptA_{\sig(i)\mu(i)}| \CP \ .
\end{equation}

Assume $\sig_o \in S_n$ and $\mu_o \in M_n$ achieve the
$\nu$-value of $|\tB|$. In case $\sig_o$ is the identity, by
Equation \Ref{proof:1}, $b_{ii} = |\tA|$,  for each $i$,
 and thus,
\begin{equation}\label{proof:4}
\bigodot_i b_{i\sig_o(i)} = \bigodot_i b_{ii} = |\tA|^n.
\end{equation}
Otherwise, when $\sig_o$ is not the identity, we write
\begin{equation}\label{proof:5}
  c := \bigodot_i    \OP \pa_{i\mu_o(i)} |\ptA_{\sig_o(i)\mu_o(i)}| \CP ,
\end{equation}
for the product that reaches the $\nu$-value of \Ref{proof:3} and
prove it always $\prec |\tA|^n$.

\parSpc \noindent  \bfem{Case I:} Assume $\mu_o \in S_n$ is a
permutation, then Formula \Ref{proof:5} can be reordered to the
form
\begin{equation}\renewcommand{\theequation}{*}\addtocounter{equation}{-1}\label{eq:str.1}
\pa_{1\mu_o(1)} |\ptA_{1\mu_o(1)}|  \ \cdots \ \pa_{n\mu_o(n)}
|\ptA_{n\mu_o(n)}|.
\end{equation}
If $(*) \succ |\tA|^n$, then it must have at least one component
$\pa_{j\mu_o(n)} |\ptA_{j\mu_o(n)}| \succ |\tA|$, but this
contradicts the maximality of $|\tA|$. On the other hand, if all
$\pa_{j\mu_o(n)} |\ptA_{j\mu_o(n)}| = |\tA|$ we get a
contradiction to the regularity of $|\tA|$. Therefore,  $(*) \prec
|\tA|^n$.

\parSpc
\noindent \bfem{Case II:} Assume $\mu_o \in \maps_n \setminus
S_n$,
  then there exist at least  two indices $i_1 \neq i_2$ for which
$ \mu_o(i_1)= \mu_o(i_2) = j_o$.  We show the existence of a
permutation $\mu_l \in S_n$ that reaches the same $\nu$-value for
Formula \Ref{proof:5} as $\mu_o$ reaches, the proof is then
completed by Case I,  applied to $\mu_l$.

For the two components $\pa_{i_1j_o}  |\ptA_{\sig_o(i_1)j_o}|$ and
$\pa_{i_2j_o}  |\ptA_{\sig_o(i_2)j_o}|$ of (\ref{proof:5}),
indexed by $ \mu_o(i_1)= \mu_o(i_2) = j_o$, we have the following
combinatorial layouts:
$$
\begin{array}{|cc|c|ccc|c|c|}
\hline
     &  & j_o  &  &  &  & j_h & \\
     &  & \bbx &  &  &  &  & \\
     &  & \bbx &  &  &  &  & \\\hline
   i_1 &  & * &  &  &  &  &\\\hline
     &  & \bbx &  &  &  &  & \\
     &  & \bbx &  &  &  &  & \\\hline
   \sig_o(i_1)  &  \bbx & \bbx & \bbx & \bbx & \bbx & \bbx & \bbx \\\hline
     &  & \bbx &  &  &  &  & \\ \hline
\end{array} \qquad
\begin{array}{|cc|c|ccc|c|c|}
\hline
     &  & j_o  &  &  &  & j_h & \\
     &  & \bbx &  &  &  &  & \\\hline
   \sig_o(i_2)  &  \bbx & \bbx & \bbx & \bbx & \bbx & \bbx & \bbx \\\hline
     &  & \bbx &  &  &  &  & \\
     &  & \bbx &  &  &  &  & \\
     &  & \bbx &  &  &  &  & \\ \hline
     i_2 &  & * &  &  &  &  &\\\hline
     &  & \bbx &  &  &  &  & \\ \hline
\end{array}
$$
The diagrams are useful to understand the modification of $\mu_o$.

Since $\mu \in M_n \setminus S_n$, there exists at least one index
$j_h \neq j_o$ in $N \setminus Im(\mu_o)$. Therefore, the
corresponding component, $\pa_{i j_h}|\ptA_{\sig_o(i) j_h}|$, is
absent in \Ref{proof:5}. Without loss of generality, we take
$\pa_{i_2 j_o} |\ptA_{\sig_o(i_2)j_o}|$ and modify it. Clearly,
$|\ptA_{\sig_o(i_2)j_o}|$ involves an entry $\pa_{\dt j_h}$, let
$i_h$ be the index for which $\sig(i_h) = j_h$. Then
$|\ptA_{\sig_o(i_2) j_o}| = \pa_{i_h j_h}
|\ptA_{\sig_o(i_2)i_h,j_oj_h}|$,
 and hence, by the maximality of $\mu_o$,
$$\pa_{i_2j_o}  |\ptA_{\sig_o(i_2)j_o}|
 = \pa_{i_2j_o}  \pa_{i_h j_h}
|\ptA_{\sig_o(i_2)i_h,j_oj_h}| = \pa_{i_h j_h} |\ptA_{i_hj_h}| \
.$$
Namely, we have specified another map $\mu_1 \in M_n$ with
$\mu_1(i_1)= j_o$, $\mu_1(i_2)= j_h$, and $\mu_1(i)= \mu_o(i)$ for
all $i \neq i_1,i_2 $. Therefore, we reduced the number of indices
sharing a same image  in $\mu_o$ to have
 $Im(\mu_o) \subset Im(\mu_1) \subseteq N$. Proceeding inductively
 we get a chain  $$Im(\mu_o) \subset Im(\mu_1)\subset \cdots \subset Im(\mu_l)=
 N,$$  the left equality is due to the finiteness  of $\maps_n$.
 Thus $\mu_l \in S_n$; the proof of Case II is then completed by Case
 I.

 So, we have showed that the identity $\sig_o$ is the single permutation
 that maximizes \Ref{proof:3}, and for which we have $\tB = \tA \itA = |\tA|\gtI
 $. Since $|\tA| \in \Real$ and $\gtI$ is regular, so is $\tB$.
 This completes  the proof of Theorem \ref{thm:inverseMatrix}
 on pseudo invertibility of matrices over $(\Trop,\TrS,\TrP)$.
\end{proof}

We push the result of  Theorem \ref{thm:inverseMatrix} further:
\begin{theorem}  For each
regular matrix $\tA \in M_{n}(\Trop)$, the products $\tA\itA$ and
$\itA\tA$ are idempotents.
\end{theorem}
\begin{proof} Writing  $\gtI = \tA\itA$, with $\gtI = (\io_{ij})$, we
prove that $\gtI = \gtI^2$.  Recall that $\ia_{ij} =
\frac{|\tA_{ji}|}{|\tA|}$, then
\begin{equation}\label{eq:multToIdim.0}
\io_{ij} = \bigoplus_k \aaa_{ik} \ia_{kj} =  \bigoplus_k
\aaa_{ik}\frac{|\tA_{jk}|}{|\tA|} =
\aaa_{ik_e}\frac{|\tA_{jk_e}|}{|\tA|} \ ,
\end{equation}
for some fixed $k_e$.  Suppose $(\gtI)^2 = (\io^{(2)}_{ij})$, then
$$
\io^{(2)}_{ij} = \bigoplus_{h} \io_{ih}\io_{hj} = \bigoplus_{h}
\OP \bigoplus_k \aaa_{ik}\frac{|\tA_{hk}|}{|\tA|} \CP \OP
\bigoplus_l \aaa_{hl}\frac{|\tA_{jl}|}{|\tA|}\CP =$$
\begin{equation}\label{eq:multToIdim.1}
 \bigoplus_{h}
 \bigoplus_k \bigoplus_l \aaa_{ik}\frac{|\tA_{hk}|}{|\tA|}
\aaa_{hl}\frac{|\tA_{jl}|}{|\tA|} \ ,
\end{equation}
and we need to prove the equality
\begin{equation}\label{eq:multToIdim.2}
 |\tA|\bigoplus_k \aaa_{ik}|\tA_{jk}| =
\bigoplus_{h} \bigoplus_k \bigoplus_l
\underbrace{\aaa_{ik}|\tA_{hk}|}_{(I)}
\underbrace{\aaa_{hl}|\tA_{jl}|}_{(II)} \ . \end{equation}

To see that $\io^{(2)}_{ij}\succeq \io_{ij}$, take $h = j$ to have
$ (II) = \bigoplus_l  \aaa_{hl}|\tA_{jl}| =  \bigoplus_l
\aaa_{jl}|\tA_{jl}| = |\tA|$.
By the way of contradiction, assume $\io^{(2)}_{ij}\succ
\io_{ij}$, and suppose $k_o$, $h_o$, and $l_o$ are the indices
reaching the $\nu$-value of $\io^{(2)}_{ij}$ in Formula
\Ref{eq:multToIdim.1}, then
$$ \aaa_{ik_o}|\tA_{h_ok_o}| \aaa_{h_ol_o}|\tA_{jl_o}| \succ
\aaa_{ik_e}|\tA_{jk_e}| |\tA| \ . $$
Clearly $\aaa_{ik_e}|\tA_{jk_e}| \succeq \aaa_{ik_o}|\tA_{jk_o}|,$
since otherwise we would have a contradiction to the maximality of
\Ref{eq:multToIdim.0}. Thus,
$$
\aaa_{ik_o}|\tA_{h_ok_o}| \aaa_{h_ol_o}|\tA_{jl_o}| \succ
\aaa_{ik_o}|\tA_{jk_o}| |\tA| \ ,
$$
and hence
\begin{equation}\label{eq:multToIdim.5}
 |\tA_{h_ok_o}| \aaa_{h_ol_o}|\tA_{jl_o}|
\succ |\tA_{jk_o}||\tA| \ .
\end{equation}
Due to the maximality of $|\tA|$, we also have $|\tA| \succeq
 \aaa_{h_o k_o}|\tA_{h_o k_o}| $, and by \Ref{eq:multToIdim.5} we
get
$$ |\tA_{h_o k_o}| \aaa_{h_o l_o}|\tA_{j l_o}| \succ
|\tA_{j k_o}||\tA| \succeq |\tA_{j k_o}| \aaa_{h_o k_o}|\tA_{h_o
k_o}| \ .$$
Namely,
\begin{equation}\label{eq:multToIdim.6}
 \aaa_{h_ol_o}|\tA_{jl_o}| \succ
 \aaa_{h_ok_o}|\tA_{jk_o}| \ .
\end{equation}
This contradicts the specification of $k_o$ as the index that
reaches the maximum for \Ref{eq:multToIdim.2}. This completes the
proof that $(\tA\itA)^2 = \tA\itA$; the case of multiplication on
left is proved in the same way.
\end{proof}

\begin{corollary}
A matrix $\tA$ is E-dense in $M_n(\Trop)$, with respect to
$U_n^{\idem}(\Trop)$, if and only if is tropically regular.
\end{corollary}

\begin{example}
 Take the regular matrix
  $$ \ptA = \vMat{1}{-1}{2}{2}, \quad  \text{ then }  \quad
     \itA = \vMat{2}{-1}{2}{1} (-3) \ ,$$ where  $|\tA| = 3$.
     (Recall that, in tropical sense, multiplying by $(-3)$ means dividing by $3$.)
     The product  $\ptA  \itA$
 is then

    $$ \vMat{1}{-1}{2}{2}
       \vMat{2}{-1}{2}{1} (-3) =
       \vMat{3}{\uuu{0}}{\uuu{4}}{3}(-3) =
       \vMat{0}{\uuu{(-3)}}{\uuu{1}}{0}  \in U^\idem_n(\Trop) \ .
       $$
On the other hand, if  we take the singular matrix
  $$ \ptA = \vMat{1}{-1}{4}{2}, \quad  \text{ then } \quad
     \itA = \vMat{2}{-1}{4}{1} \uuu{(-3)},$$ where here  $|\tA| = \uuu{3}$.
   Computing the product $\ptA  \itA$ we get
    $$
        \vMat{1}{-1}{4}{2}
        \vMat{2}{-1}{4}{1} \uuu{(-3)} =
        \vMat{\uuu{3}}{\uuu{0}}{\uuu{6}}{\uuu{3}} \uuu{(-3)} \notin U_n(\Trop) \ ,$$
        which is not a regular matrix,
        and therefore $\ptA  \itA \notin U_n(\Trop)$.
\end{example}

 A few immediate conclusions are derived from our  last
 results:
\begin{corollary}\label{thm:detOfinverseMatirx}
Assume $\tA$ is a regular matrix, then
\begin{enumerate}
    \item $Adj(\tA)$ is also regular; \pSkip
    \item $|\tA| = (|\itA|)^{-1}$, and if $\tA =  \itA$ then $|\tA| = |\itA| =
    0$. \pSkip
\end{enumerate}
\end{corollary}
\begin{proof} The first assertion is obvious.
 $\tA$, $\itA$, and $\tA \itA$ are all regular, then by Theorem \ref{thm:detofMult}
 $ |\tA|  |\itA| = |\gtI| = 0$ and hence $|\tA| =
|\itA|^{-1}$.
\end{proof}
\noindent The converse assertion of $(2)$ is not true; for
example, take the matrix
$$\tA  = \vMat{-1}{-2}{-2}{1}, \quad \text{ then } \quad
  \itA = \vMat{1}{-2}{-2}{-1} \ .$$
Although $|\tA| = |\itA| =0$, we have  $\tA \neq \itA$.

\begin{remark} Contrary to the classical theory of matrices over fields,
tropically,  the relation $\invA{(\tA \tB)} = \itB \itA$ does not
hold true; for example, take the regular matrix as in
\Ref{exmp:DetMultiplication}, then
\begin{equation*}
   \tA = \vMat{1}{1}{2}{3}, \quad  \itA =
   \vMat{3}{1}{2}{1}(-4), \quad \text{and} \quad  \itA \itA =
                         \vMat{6}{4}{5}{3}(-8).
 \end{equation*}
On the other hand, $\tA^2$ is not regular, cf. Remark
\ref{obs:regProp}, and the computation of $Adj(\tA^2)/|\tA^2|$
yields
\begin{equation*}
  \invA{(\tA \tA)} =
  \vMat{6}{4}{5}{3}\uuu{(-9)}, \quad \text{where} \quad  \tA^2 =
  \vMat{3}{4}{5}{6} \ ;
\end{equation*}
this shows that $(\itA)^2 \neq \invA{(\tA \tA)} $.
\end{remark}

 \subsection{Matrices with real entries}\label{sec:Observation}

 Denoting by $M_n(\eReal)$ the semiring of matrices over $(\eReal, \max,+ \
 )$, the epimorphism $\epiToMaxPlus :  (\Trop,\TrS,\TrP) \to (\eReal, \max,+ \
 ) $, cf.~\Ref{eq:homomorphismsMax2Extended},  induces in the standard way the
 epimorphism $$\epiToMaxPlus_*: M_n(\Trop) \  \To \ M_n(\eReal)$$
 of matrix semirings. We write  $\epiToMaxPlus_*(\tA)$ for the
 image of $\tA \in  M_n(\Trop) $ in $ M_n(\eReal)$.  Conversely, set-theoretic, $M_n(\eReal) \subset
 M_n(\Trop)$.

\begin{proposition}\label{thm:exactRelationOfInverse}
Suppose $\tA \in  M_n(\Trop)$ is  regular, where both $\tA$ and
$\itA$ have only real entries, $\tA  \itA = \gtI'$, and $
   \itA  \tA = \gtI''$. Then
$$
  \epiToMaxPlus_*(\gtI'  \ptA) = \tA, \ \ \epiToMaxPlus_*(\itA   \gtI') = \itA, \ \
  \epiToMaxPlus_*(\gtI''  \itA) = \itA,  \ \text{and } \ \epiToMaxPlus_*(\ptA   \gtI'') = \ptA.
$$
\end{proposition}
\begin{proof} We prove the relation $\epiToMaxPlus_*(\gtI'  \ptA) =
\tA$. Letting $\gtI' = (\io_{ij})$, we  show that
$$\epiToMaxPlus ( \bigoplus_k \io_{ik}  \aaa_{kj}) = \aaa_{ij}
,$$
for all $i,j$. Recall that $
  \io_{ik} =
   \OP \bigoplus_{h} \pa_{ih}  |\ptA_{jh}| \CP |\tA|^{-1}
$, cf. Formula~\Ref{proof:2}, and $\io_{ik} \in \etUnit$ whenever
$i \neq k$. Composing together, we get
\begin{equation}\renewcommand{\theequation}{*}\addtocounter{equation}{-1}\label{eq:str.1}
\bigoplus_k \OP \bigoplus_{h} \pa_{ih} |\ptA_{jh}||\tA|^{-1} \CP
\aaa_{kj} =  \OP \bigoplus_{k,h} \pa_{ih}  |\ptA_{jh}| \aaa_{kj}
\CP  |\tA|^{-1} \  .
\end{equation}
Using  Formulas \Ref{proof:4} and \Ref{proof:5}, we see that the
maximal value of $|\ptA_{jh}| \aaa_{kj}$ is attained when $k = h =
j$ and it is $|\tA|$. Thus,
$$\epiToMaxPlus\OP(*)\CP =
\epiToMaxPlus(\pa_{ij}  |\ptA_{jj}|  \aaa_{jj} |\tA|^{-1})
=\pa_{ij} \ .
$$
The other relations are proved in the same way.
\end{proof}

\begin{remark}
In the sense of Proposition \ref{thm:exactRelationOfInverse}, the
matrices $\gtI'$ and $\gtI''$ are pseudo right/left identities of
$\tA$ and $\itA$ respectively.
\end{remark}

Pushing the results of Proposition
\ref{thm:exactRelationOfInverse} forward, we conclude:
 \begin{corollary}\label{thm:exactRelationOfInverse}
Suppose $\tA \in  M_n(\Trop)$ is  regular.  Let $\tA  \itA =
\gtI'$ and $
   \itA  \tA = \gtI''$;
then
$$
  \epiToMaxPlus_*(\gtI'  \ptA) = \epiToMaxPlus_*(\tA), \ \
  \epiToMaxPlus_*(\itA   \gtI') = \epiToMaxPlus_*(\itA), \ \
  \epiToMaxPlus_*(\gtI''  \itA) = \epiToMaxPlus_*(\itA),  \
   \text{and } \ \epiToMaxPlus_*(\ptA   \gtI'') = \epiToMaxPlus_*(\ptA).
$$
\end{corollary}

\begin{example}
Let
  $$ \ptA = \vMat{1}{1}{2}{3} \, \quad  \text{ then } \quad
     \itA = \vMat{-1}{-3}{-2}{-3} \quad \text{and} \quad
   \ptA \itA  = \gtI' =
     \vMat{0}{\uuu{(-2)}}{\uuu{1}}{0}.$$
Computing the products we have
  $$
\begin{array}{lllll}
   \gtI'  \ptA  & = & \vMat{0}{\uuu{(-2)}}{\uuu{1}}{0}
   \vMat{1}{1}{2}{3} & = & \vMat{1}{\uuu{1}}{\uuu{2}}{3}, \\ [4mm]
 \itA  \gtI' & = & \vMat{-1}{-3}{-2}{-3}
    \vMat{0}{\uuu{(-2)}}{\uuu{1}}{0} & =&
    \vMat{-1}{\uuu{(-3)}}{\uuu{(-2)}}{-3} \ ,
  \end{array}
  $$
  and it easily verify that  $\epiToMaxPlus_*( \gtI'  \ptA) = \ptA$ and $\epiToMaxPlus_*(\itA  \gtI') =
  \itA$.
\end{example}


\end{document}